\documentstyle[12pt,psfig]{article}
\psfigurepath{ps}

\newfont{\ghot}{eufm10 scaled \magstep 2} 
\newfont{\sets}{msbm10 scaled \magstep 2} 
\def\bl{\rule[-1mm]{2.4mm}{2.4mm}}
\newtheorem{thrm}{T~h~e~o~r~e~m}
\newtheorem{lmm}{L~e~m~m~a}

\begin{document}
\title{Integral equations PS-3 and moduli of pants.}
\author{\copyright 2006 ~~~~A.B.~Bogatyr\"ev
\thanks{Supported by RFBR grant 05-01-01027 and grant MD4798.2007.01}}
\date{}
\maketitle

We introduce the  new type of constructive
\index{representation!pictorial} {\it pictorial} representations
for the solutions of the following spectral singular
Poincar\'e--Steklov (PS for brevity) integral equation
\index{equation!singular~integral}
\index{equation!Poincar\'e-Steklov}
\begin{equation}
\label{PSE}
\lambda~
V.p.\int\limits_I
\frac{u(t)}{t-x}\,dt -
V.p.\int\limits_I
\frac{u(t)\,dR(t)}{R(t)-R(x)}
=const,\quad x\in I:=(-1,~1),
\end{equation}
where $\lambda$ is the spectral parameter; $u(t)$ is the unknown
function; $const$ is independent of $x$. The functional parameter
$R(t)$ of the equation is a smooth {\it  nondegenerate} change of
variable on the interval $I$:
\begin{equation}
\frac d{dt}
R(t)\neq 0,
\qquad t\in [-1,1].
\label{ND}
\end{equation}

\section{Introduction}
H.Poincar\'{e} (1896) and V.A.Steklov (1901) were the first who placed the spectral parameter
to the boundary conditions of the problem for an elliptic operator.
Later it became a popular technique for the analysis and optimization in diffraction problems
\cite{AKSV}, (thermo)conductivity of composite materials,  simple 2D model of oil extraction
 etc.

\subsection{Spectral Boundary Value Problem}
\index{spectral!parameter}
\index{spectral!problem}
Let a domain in the plane be subdivided into two simply connected
domains $\Omega_1$ and $\Omega_2$ by a smooth simple arc $\Gamma$.
We are looking for the values of the spectral parameter $\lambda$
when the following problem has a nonzero solution:
\begin{equation}
\begin{picture}(30,30)
\put(0,10){\oval(40,30)}
\put(0,-5){\oval(10,30)[tl]}
\put(0,25){\oval(10,30)[br]}
\put(10,3){$\Omega_1$}
\put(-12,13){$\Omega_2$}
\put(6,16){$\Gamma$}
\end{picture}
\begin{array}{rcllrcl}
\Delta U_1\!&=&\!0 & {\rm in} \ \Omega_1;
& U_1\!&=&\!0 ~{\rm on} \ \partial\Omega_1\setminus\Gamma; \\
\Delta U_2\!&=&\!0 & {\rm in} \ \Omega_2;
& U_2\!&=&\!0 ~{\rm on} \ \partial\Omega_2\setminus\Gamma;\\
{~~}U_1\!&=&\!U_2  \ & {\rm on} \ \Gamma;  &&&\\
-\lambda\displaystyle{\frac{\partial U_1}{\partial n}}\!&=&
\!\displaystyle{\frac{\partial U_2}{\partial n}} &{\rm
on} ~\Gamma,  &&&\\
\end{array}
\label{BVP}
\end{equation}

\index{Boundary value problem} \index{Domain decomposition method}
Spectral problems of this type naturally arise e.g. in the
justification and optimization of a {\it domain decomposition
method} for the solution of a boundary value problem for Laplace
equation. It is easy to show that the   eigenfunctions and the
eigenvalues of the problem (\ref{BVP}) are correspondingly the
critical points and critical values of the functional (the so
called generalized Rayleigh ratio) \index{eigenvalue}
\index{eigenfunction} \index{generalized Rayleigh ratio}
\begin{equation}
F(U)=\frac{\int_{\Omega_2}|grad ~U_2|^2 ~d\Omega_2}
{\int_{\Omega_1}|grad ~U_1|^2 ~d\Omega_1},
 \qquad U\in H_{oo}^{1/2}(\Gamma),
\label{Rayleigh}
\end{equation}
  where $U_s$ is the harmonic continuation of the function $U$
  from the interface $\Gamma$ to the domain $\Omega_s$, $s=1,2$,
  vanishing at the outer boundary of the domain.

  The boundary value problem (\ref{BVP}) is equivalent to a
  certain Poincar\'{e}--Steklov equation. Indeed, let $V_s$ be the harmonic
  function conjugate to $U_s$, $s=1,2$. From the Cauchy--Riemann \index{equations!Cauchy-Riemann}
  equations and the relations on $\Gamma$ it follows that the tangent
  to the interface derivatives of $V_1$ and $V_2$ differ by the
  same factor $-\lambda$. Integrating along $\Gamma$ we get
\begin{equation}
\lambda V_1(y)+V_2(y)=const,
\quad y\in\Gamma.
\label{HarmConj}
\end{equation}
    The boundary values of conjugate functions  harmonic  in the
  half-plane are related by a Hilbert transformation. To reduce our
  \index{Hilbert transformation}
  case to this model we consider a conformal mapping
  $\omega_s(y)$ from $\Omega_s$ to the open upper half-plane
  $\mbox{\sets H}$ with normalization $\omega_s(\Gamma)=I$,
  $s=1,2$. Now equation (\ref{HarmConj}) may be rewritten as
\index{conformal mapping}

$$
-\frac{\lambda}{\pi}~
V.p.\int\limits_I
\frac{U_1(\omega_1^{-1}(t))}{t-\omega_1(y)}\,dt -
\frac1{\pi}~
V.p.\int\limits_I
\frac{U_2(\omega_2^{-1}(t'))}{t'-\omega_2(y)}\,dt'
=const,
\quad y\in\Gamma.
$$
Introducing the new notation $x:=\omega_1(y)\in I$;
$R:=\omega_2\circ\omega_1^{-1}:~~I\to\Gamma\to I$;
$u(t):=U_1(\omega_1^{-1}(t))$ and the change of variable $t'=R(t)$
in the second integral we arrive at the Poincar\'{e}--Steklov equation
(\ref{PSE}). Note that in this context $R(t)$ is the decreasing
function on $I$.

\subsection{Some Known Results}

  The natural way to study integral equations is  operator
  analysis. This discipline allows to obtain
  for the {\it  smooth nondegenerate} change of variables $R(x)$
 the following results \cite{Bog0}:

\begin{itemize}
\item The spectrum is discrete; the eigenvalues are positive and converge to $\lambda=1$.
\item  $\sum_{\lambda\in Sp}|\lambda-1|^2<\infty$ (a constructive estimate in terms of $R(x)$ is given)
\item The eigenfunctions $u(x)$ form a basis in the Sobolev space $H_{oo}^{1/2}(I)$.
\end{itemize}

\subsection{Goal and Philosophy of the Research}
The approach of complex geometry for the same integral equation gives
different types of results.
For quadratic $R(x)=x+(2C)^{-1}(x^2-1)$, $C>1$, the eigenpairs were found explicitly \cite{Bog1}:
$$
u_n(x)=\sin
\left[
\frac{n\pi}{K'}\int\limits_1^{(C+x)/(C-1)}(s^2-1)^{-1/2}
(1-k^2s^2)^{-1/2}ds
\right],
$$
$$
\lambda_n=1+1/\cosh 2\pi\tau n,
 \quad n=1,2\dots,
$$
where $\tau=K/K'$ is the ratio of the complete elliptic integrals of modulus $k=(C-1)/(C+1)$.
  \index{integral!elliptic}
  Now we are going to give {\it  constructive representations} for
  the eigenpairs $\{\lambda,u(x)\}$ of the integral equation with
  $R(x)=R_3(x)$ being a rational function of degree 3. Equation
  (\ref{PSE}) itself will be called PS-3 in this case.

  The notion of a constructive representation for the solution
  should be however specified. Usually this means that we
  restrict the search for the solution to a certain class of
  functions such as rational, elementary, abelian, quadratures,
  the Umemura classical functions, etc. The history of mathematics
  knows many disappointing results when the solution of the
  prescribed form does not exist. Say, the diagonal of the square
  is not commensurable with its side, generic algebraic equations
  cannot be solved in radicals, linear ordinary differential
  equations usually cannot be solved by quadratures, Painlev\'e
  equations cannot be solved by Umemura functions. Nature
  \index{classical Umemura function} \index{equation!Painleve}
  always forces us to introduce new types of transcendent
  objects to enlarge the scope of search. The study of new
  transcendental functions constitutes the progress of mathematics. This
  research philosophy goes back to H.~Poincar\'{e} \cite{Poi}. From the
  philosophical point of view our goal is to  disclose the nature
  of emerging transcendental functions in the case of PS-3 integral
  equations.

  \subsection{Brief Description of the Result}

  The rational function $R_3(x)$ of  degree three is explicitly related
  to a {\it  pair of pants} in section \ref{PantsOfR3}. On the
  other hand, given a spectral parameter $\lambda$ and two
  auxiliary real parameters, we explicitly construct in section
  \ref{Result} another pair of pants which additionally depend on
  \index{pants!pair of} \index{pants!moduli}
  two integers. When the above two pants are
  conformally equivalent, $\lambda$ is the eigenvalue of the
  integral equation PS-3 with parameter $R_3(x)$. Essentially,
  this means that to find the spectrum of the given integral
  equation (\ref{PSE}) one has to solve three transcendental
  equations involving three {\it  moduli of pants}.

  Whether this representation of the solutions may be considered
  as constructive or not is a matter of discussion. Our approach
  to the notion of a constructive representation is utilitarian:
  the more we learn about the solution from the given
  representation the more constructive is the latter.  At least
  we are able to obtain valuable features of the solution: to
  determine the number of zeroes of the eigenfunction $u(t)$, to find the
  exact locus for the spectra and to show the discrete
  mechanism of generating the eigenvalues.

\section{Description of the Main Result}

  The shape of the two domains $\Omega_1$ and $\Omega_2$ defines the
  variable change $R(x)$ only up to a certain two-parametric
  deformation. One can easily check that the {\it  gauge
  transformation} $R\to L_2\circ R\circ L_1$, where the linear
  fractional function $L_s(x)$ keeps the segment $[-1,1]$, does
  not affect the spectrum of equation (\ref{PSE}) and induces
  only the change of the argument for its eigenfunctions:
  $u(x)\to$ $u\circ L_1(x)$. For this reason we do not distinguish
  between two PS equations with their functional parameters $R(x)$
  related by the gauge transformation. \index{gauge transformation}

  The space of equivalence classes of equations PS-3
  has real dimension $3=7-2-2$ and several components with different
  topology of the functional parameter $R_3$. In the present
  paper we study for brevity only one of the components, the choice is
  specified in section \ref{ComponentDef}.

\subsection{Topology of the Branched Covering}
\label{R3Topology}
  In what follows we consider {\it  rational degree three}  functions
  $R_3(x)$ with {\it  separate real critical values} different from
  $\pm1$.  The rational function $R_3(x)$ defines a  3-sheeted branched covering of a Riemann sphere by another Riemann
  sphere. The Riemann--Hurwitz formula suggests that $R_3(x)$
  has four separate branch points $a_s$, $s=1,\ldots,4$. This
  means that every value $a_s$ is covered by a critical (double) point
  $b_s$, and an ordinary point $c_s$.
  \index{critical point}\index{critical value}\index{Riemann-Hurwitz formula}

  Every point $y\neq a_s$ of the extended real axis
  $\hat{\mbox{\sets R}}:=\mbox{\sets R}\cup\{\infty\}$ belongs to
  exactly  one of two types. For the type (3:0) the pre-image
  $R_3^{-1}(y)$ consists of three distinct real points. For the
  type (1:2) the pre-image consists of a real and two complex
  conjugate points. The type of the point remains locally constant on
  the extended real axis and changes when we step over the branch
  point. Let the branch points $a_s$ be enumerated in the
  natural cyclic order of $\hat{\mbox{\sets R}}$ so that the
  intervals $(a_1,a_2)$ and $(a_3,a_4)$ are filled with the
  points of the type (1:2). Later we will specify the way to
  exclude the relabeling $a_1\leftrightarrow a_3$,
  $a_2\leftrightarrow a_4$ of branch points.


\begin{figure}[ht!]
\begin{picture}(165,30)
\put(110,0){
\begin{picture}(70,20)

\thicklines
\put(10,10){\circle*{1}}
\put(20,10){\circle*{1}}
\put(10,10){\line(1,0){10}}
\put(8,12){$a_1$}
\put(18,12){$a_2$}
\put(30,10){\circle*{1}}
\put(40,10){\circle*{1}}
\put(30,10){\line(1,0){10}}
\put(28,12){$a_3$}
\put(38,12){$a_4$}

\thinlines
\multiput(0,10)(2,0){25}{\line(1,0){1}}
\put(51,9){$\hat{\mbox{\sets R}}$}
\put(10,25){\circle{4}}
\put(9,24.2){$y$}

\end{picture}
}
\put(5,0){
\begin{picture}(70,20)
\thicklines
\put(10,10){\circle*{1}}
\put(20,10){\circle*{1}}
\put(10,10){\line(1,0){10}}
\put(8,12){$c_4$}
\put(18,12){$c_3$}
\put(50,10){\circle*{1}}
\put(60,10){\circle*{1}}
\put(50,10){\line(1,0){10}}
\put(48,12){$c_2$}
\put(58,12){$c_1$}
\put(15,10){\oval(30,20)}
\put(0,10){\circle*{1}}
\put(30,10){\circle*{1}}
\put(1,12){$b_1$}
\put(26,12){$b_2$}
\put(55,10){\oval(30,20)}
\put(40,10){\circle*{1}}
\put(70,10){\circle*{1}}
\put(41,12){$b_3$}
\put(66,12){$b_4$}
\thinlines
\multiput(-5,10)(2,0){40}{\line(1,0){1}}
\put(76,9){$\hat{\mbox{\sets R}}$}
\put(90,10){$\longmapsto$}
\put(90,15){$R_3(x)$}
\put(75,25){\circle{4}}
\put(74,24){$x$}

\end{picture}}
\end{picture}
\caption[]
{\normalsize The topology of the covering $R_3$ with real branch points}
\label{CoverTopology}
\end{figure}
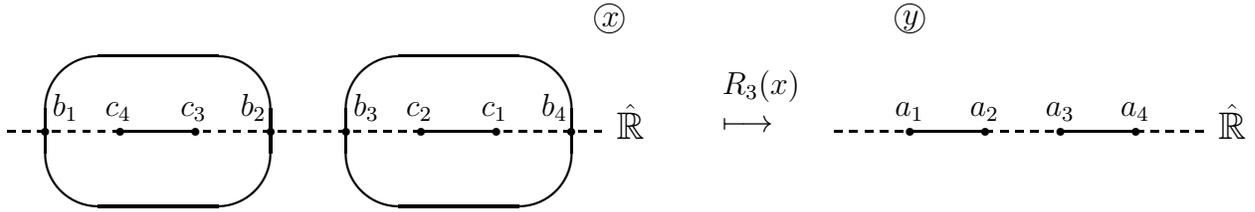

  The total pre-image $R_3^{-1}(\hat{\mbox{\sets R}})$ consists of
  the extended real axis  and two pairs of complex conjugate arcs
  intersecting $\hat{\mbox{\sets R}}$ at the points $b_1$, $b_2$,
  $b_3$, $b_4$ as shown at the left side of Fig.
  \ref{CoverTopology}. The compliment of this pre-image on the
  Riemann sphere has six components, each of them is mapped 1-1
  onto the upper or lower half-plane.

  \subsubsection{The Component in the Space of Equations}
  \label{ComponentDef}
  The  nondegeneracy condition (\ref{ND})
  forbids that any of critical points $b_s$ be inside the segment of
  integration $[-1,1]$. In what follows we consider the case when
  the latter segment lies in the annulus bounded by two ovals
  passing through the critical points $b_s$. Possibly relabeling
  the branch points we assume that $[-1,1]\subset$ $(b_2,b_3)$.

  Other components in the space of PS-3 integral equations are treated in \cite{Bog3}.

\subsection{Pair of Pants}
\label{PantsOfR3}\index{pants!pair of}\index{pants!conformal class of}
  For obvious reason {\it  a pair of pants} is the name for the
  Riemann sphere with three holes in it. Any pair of pants may
  be conformally mapped to $\hat{\mbox{\sets C}}:=$
  ${\mbox{\sets C}}\cup\{\infty\}$ with  three
  nonintersecting real slots. This mapping is unique up to the
  real linear-fractional transformation of the sphere. The conformal
  class of pants with labeled boundary components depend on three
  real parameters varying in a cell.

{\bf Definition}
~~To the variable change $R_3(x)$ we associate the pair of pants
\begin{equation}
\label{R3Pants}
{\cal P}(R_3):=Closure
\left(
\hat{\mbox{\sets C}}
\setminus
\{
[-1,1]\cup [a_1,a_2]\cup [a_3,a_4]
\}
\right).
\end{equation}
{\it  Closure} here and below is taken with respect to the
intrinsic spherical metrics when every slot acquires two sides.
 Boundary components of the pair of pants are colored (labeled) in accordance with the palette:
\begin{center}
\begin{tabular}{ll}
$[-1,1]$&     -- "red",\\
$[a_1,a_2]$& -- "blue",\\
$[a_3,a_4]$& -- "green".
\end{tabular}
\end{center}

The conformal class of pants (\ref{R3Pants}) depends only on the
equivalence class of integral equations.
To simplify the statement of our result we assume that infinity lies strictly inside the pants
(\ref{R3Pants}) which is not a loss of generality -- we can always apply a suitable gauge transformation of
the parameter $R_3(x)$.

\subsubsection{Reconstruction of $R_3(x)$ from the Pants}
  Here we show that the branched covering map $R_3(x)$ with given
  branch points $a_s$, $s=1,\dots,4$, is essentially unique.
  A possible ambiguity is due to the conformal motions
  of the covering Riemann sphere.

  Let $L_a$ be the unique linear-fractional map sending the
  critical values $a_1$, $a_2$, $a_3$, $a_4$ of $R_3$ to respectively $0$,
  $1$, $a>1$, $\infty$. The conformal motion $L_b$ of the
  covering Riemann sphere sends the critical points $b_1$, $b_2$,
  $b_3$, $b_4$ (unknown at the moment) to respectively $0$, $1$,
  $b>1$, $\infty$. The function $L_a\circ R_3\circ
  L_b^{-1}$ with normalized critical points and critical values
  takes a simple form

  $$
  \widetilde{R_3}(x)=x^2L(x)
  $$
  with a real linear fractional function $L(x)$ satisfying
  the restrictions:
$$
\begin{array}{ll}
L(1)=1,&L'(1)=-2,\\
L(b)=a/b^2,&L'(b)=-2a/b^3.
\end{array}
$$
We got four equations for three parameters of $L(x)$ and the unknown $b$.
The first two equations suggest the following expression for $L(x)$
$$
L(x)=1+2\frac{(c-1)(x-1)}{x-c}.
$$
  Another two are solved parametrically in terms of parameter $c$:
  $$
  b=c\frac{3c-3}{2c-1}; \qquad a=c\frac{(3c-3)^3}{2c-1}.
  $$
  Both functions $b(c)$ and $a(c)$ increase from $1$ to $\infty$
  when the argument $c\in (1/3,1/2)$. So, given $a>1$ we find the unique $c$
in just specified limits, and therefore the mapping $\widetilde{R_3}(x)$.
  Now we can restore the linear fractional map  $L_b$. The
  inverse image $\widetilde{R_3}^{-1}$ of the segment $L_a[-1,1]$
  consists of three disjoint segments. For {\it  our case} we
  choose the (unique) component of the pre-image belonging to the
  segment $[1,b]$. The requirement: $L_b$ {\it maps $[-1,1]$ to
  the chosen segment} determines $R_3(x)$ up to a gauge
  transformation.

\subsection{Another Pair of Pants} \label{Result} \index{annulus}
For real $\lambda\in(1,2)$ we consider an annulus $\alpha$ depending on $\lambda$ bounded by $\varepsilon\hat{\mbox{\sets R}}$,
$\varepsilon:=\exp(2\pi i/3)$, and the circle
\begin{equation}
\label{circle}
C:=
\{p\in\mbox{\sets C}:
\quad |p-\mu^{-1}|^2=\mu^{-2}-1\},
\qquad \mu:=\sqrt{\frac{3-\lambda}{2\lambda}}\in(\frac12,1).
\end{equation}
Another annulus bounded by the same circle $C$ and
$\varepsilon^2\hat{\mbox{\sets R}}$ is denoted by $\overline{\alpha}$.
Note that for the considered values of $\lambda$ the circle $C$ does not intersect the lines
$\varepsilon^{\pm1}\mbox{\sets R}$. We paint the boundaries of our annuli in the following way:

\begin{center}
\begin{tabular}{ll}
$C$&     -- "red",\\
$\varepsilon\hat{\mbox{\sets R}}$& -- "green",\\
$\varepsilon^2\hat{\mbox{\sets R}}$& -- "blue".
\end{tabular}
\end{center}

Given $\lambda$ in the specified above limits, real $h_1,~h_2$
and nonnegative integers $m_1,m_2$, we define three pairs of pants
${\cal P}_s(\lambda,h_1,h_2|m_1,m_2)$ of different fashions  $s=1,2,3$.

\noindent{\bf Fashion~1:}\\
${\cal P}_1(\lambda,h_1,h_2|m_1,m_2):=~~$ $m_1\alpha~~+~~m_2\overline{\alpha}~~+$\\[3mm]
\psfig{figure=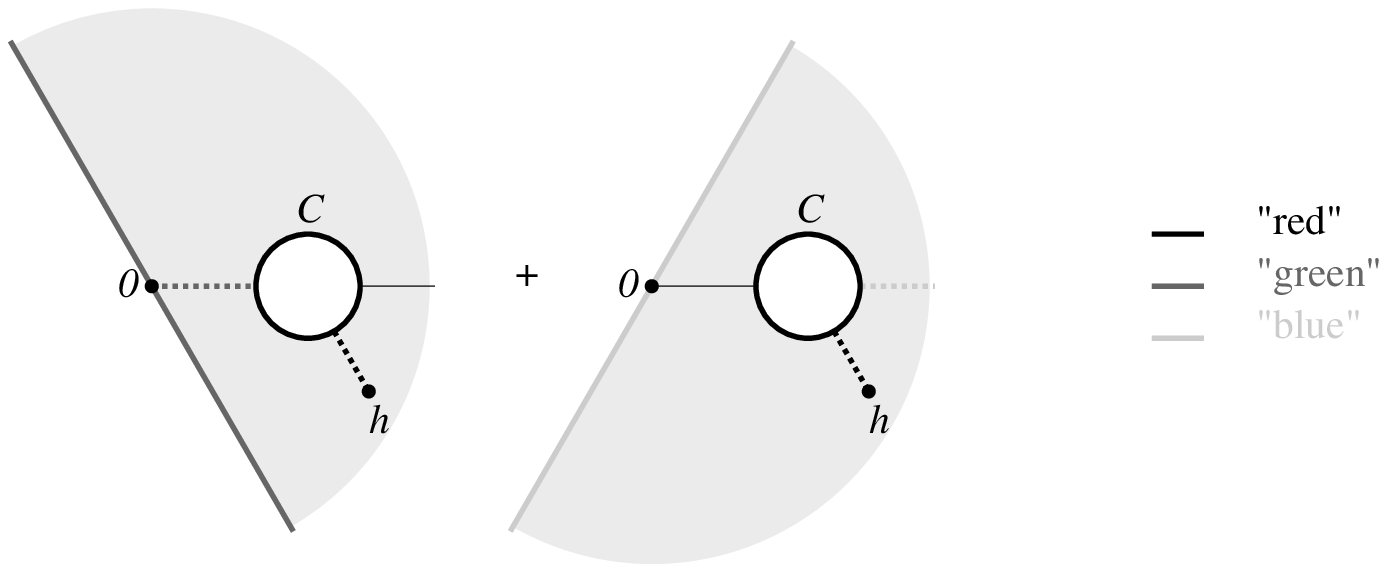}

\begin{figure}[ht!]
\centerline{\psfig{figure=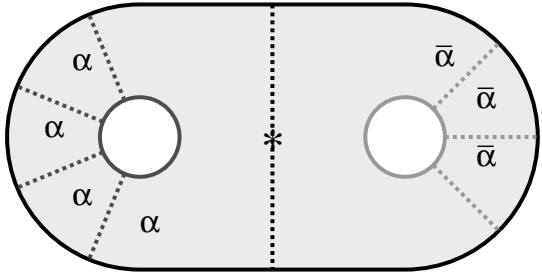}} \caption[] {\normalsize The
scheme for sewing pants ${\cal P}_1(\lambda,h_1,h_2|3,2)$.
Asterisk is the critical point of $p(y)$.} \label{Sew1}
\end{figure}

  The operations $'+'$ here stand for a certain surgery. First of all
  take two annuli $\alpha$ and  $\overline{\alpha}$ and cut them
  along the same segment (dashed red line in the figure above)
  starting at the point $h:=h_1+ih_2$ from the interior of
  $\alpha\cap\overline{\alpha}$ and ending at the circle $C$. Now
  glue the left bank of one cut to the right bank of the other.
  The resulting two sheeted surface (called
  {\it \"Uberlagerungsfl\"ache} in the following)
  \index{\"Uberlagerungsfl\"ache}
  will be the pair of pants ${\cal P}_1(\lambda,h_1,h_2|0,0)$. It
  is possible to modify the obtained surface sewing several
  annuli to it. Cut the annulus $\alpha$ contained in the pants
  and $m_1$ more copies of this annulus along the same segment
  (shown by the dashed green line in the figure above) connecting
  the boundaries of the annulus. The left bank of the cut on
  every copy of $\alpha$ is identified with the right bank of the
  cut on another copy so that all copies of the annulus are glued
  in one piece. A similar procedure may be repeated for the annulus
  $\overline{\alpha}$ (cut along the dashed blue line). The
  scheme for sewing together fashion 1 pants from the patches $\alpha$, $\bar\alpha$
  when $m_1=3$ and $m_2=2$ is shown in Fig. \ref{Sew1}.

\pagebreak

  \noindent{\bf Fashions~ 2 and 3}:\\
${\cal P}_2(\lambda,h_1,h_2|m_1,m_2):=$
\hfill
${\cal P}_3(\lambda,h_1,h_2|m_1,m_2):=$\\
$m_1\alpha~~+~~m_2\overline{\alpha}~~+~~$
\hfill
$m_1\alpha~~+~~m_2\overline{\alpha}~~+~~$\\
\psfig{figure=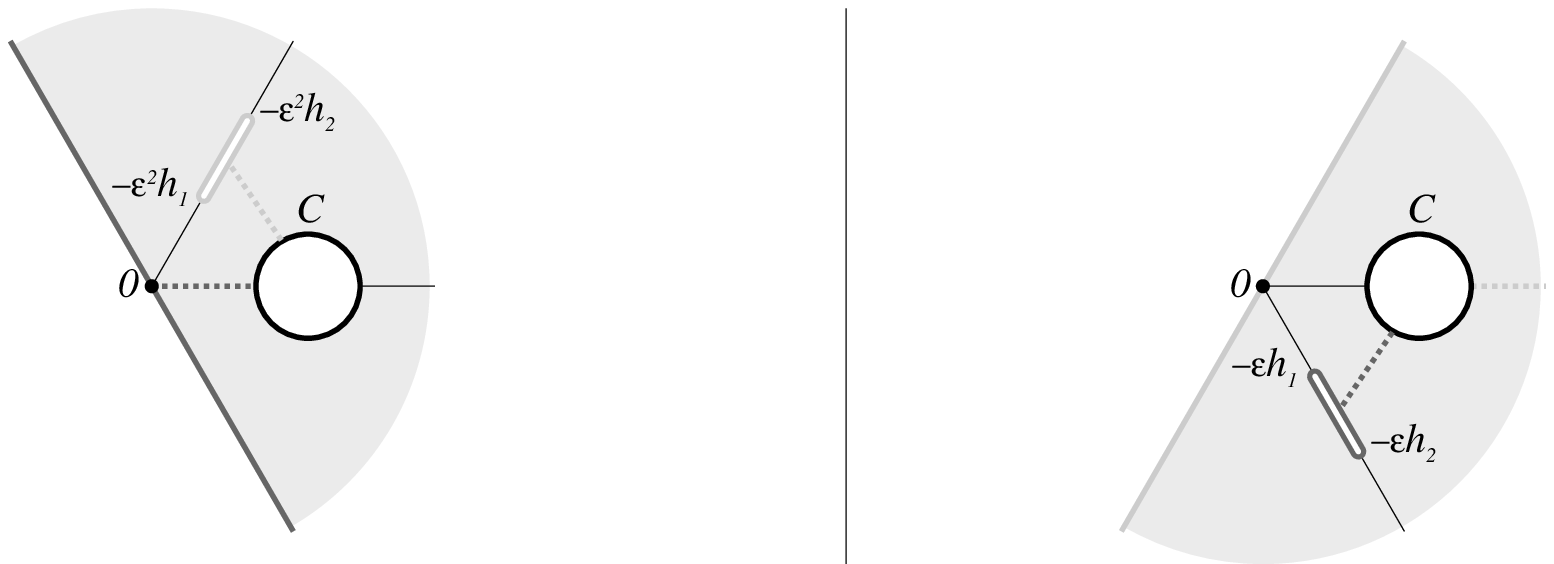}

\begin{figure}[h]
\centerline{\psfig{figure=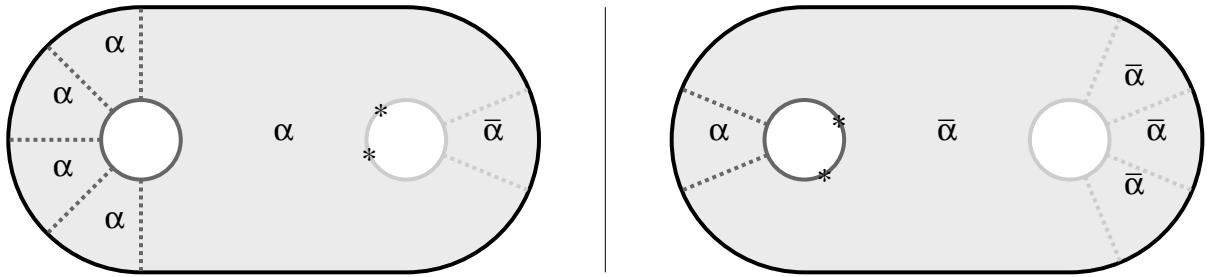}} \caption[] {\normalsize The
scheme for sewing pants ${\cal P}_2(\lambda,h_1,h_2|4,1)$ ({\bf
left}); and ${\cal P}_3(\lambda,h_1,h_2|1,3)$ ({\bf right}).
Asterisks are the critical points of the mapping $p(y)$.}
\label{Sew23}
\end{figure}

The pair of pants ${\cal P}_2(\lambda,h_1,h_2|0,0)$ (resp. ${\cal
P}_3(\lambda,h_1,h_2|0,0)$) by definition is the annulus $\alpha$
(resp. $\overline{\alpha}$) with  removed
  segment $-\varepsilon^2[h_1,h_2]$ (resp. $-\varepsilon[h_1,h_2]$),
  $0<h_1<h_2<\infty$.   As in the previous case those pants may be modified by sewing
  in several annuli $\alpha$, $\overline{\alpha}$. The scheme of cutting and gluing is shown
  in Fig. \ref{Sew23}

  \subsubsection{Remarks on the Constructed Pairs of Pants} \label{Rmk1}
   {\bf 1.} The limiting case of the first fashion of the pants when the branch point $h_1+ih_2$ tends to
  $\varepsilon^{\pm1}\mbox{\sets R}$ coincides with the limiting cases
  of the two other fashion pants when $h_1=h_2>0$:
  \begin{equation}
  \begin{array}{cc}
  {\cal P}_1(\lambda, -Re(\varepsilon^2h),-Im(\varepsilon^2h)|m_1,m_2)=&
  {\cal P}_2(\lambda, h,h|m_1,m_2+1),\\
  {\cal P}_1(\lambda, -Re(\varepsilon h),-Im(\varepsilon h)|m_1,m_2)=&
  {\cal P}_3(\lambda, h,h|m_1+1,m_2),\\
  \end{array}
  \label{InterMediate}
  \end{equation}
where parameter $h>0$.
We denote those intermediate cases as ${\cal P}_s(\lambda, h|m_1,m_2)$, $s=12, 13$ respectively.

  {\bf 2.} The surgery procedure of sewing annuli e.g. to the pants (known as {\it  grafting of projective structures})
  was designed by B.~Mas\-kit (1969), D.~Hej\-hal (1975) and
  D.~Sul\-livan-W.~Thur\-ston
  (1983), see also W.~Gold\-man (1987). \index{projective structure!grafting}

  {\bf 3.} \label{Rmk3} Every pair of pants ${\cal P}_s(\lambda,
  h_1,h_2|m_1,m_2)$ may be conformally mapped to the sphere with
  three real slots, i.e., pants of the type (\ref{R3Pants}). Let
  $p(y)$ be the inverse mapping. We observe that $p(y)$ has
  exactly one critical point in the pants ${\cal P}(R_3)$,
  counting {\it multiplicity and weight}. For the fashion $s=1$ this
  point lies  strictly inside the pants and is mapped to
  $h=h_1+ih_2$. For the case $s=2$ (resp. $s=3$) there will be
  two simple critical points of $p(y)$ on the blue (resp. green)
  \index{critical point}
  boundary component of the pants which are mapped to the points
  $-\varepsilon^2 h_1$, $-\varepsilon^2 h_2$ (resp.$-\varepsilon
  h_1$, $-\varepsilon h_2$). Finally, for the intermediate case
  (see remark 1) there will be a double critical point on the
  boundary. The {\it  multiplicity} of the critical point on the
  boundary should be calculated with respect to the local
  parameter of the {\it double} of pants ${\cal P}(R_3)$:
  \begin{equation}
  \label{M}
  M:=\{w^2=(y^2-1)\prod\limits_{s=1}^4(y-a_s)\},
  \end{equation}
  e.g., at the
  endpoint $a=\pm1, a_1,\dots$, $a_4$  of the slot this local
  parameter is $\sqrt{y-a}$. We consider the critical points on the boundary
  with the {\it weight} $\frac12$.
  \index{pants!double of}

\subsection{Main Theorem}

  Later we explain that {\it  real} eigenfunctions of the integral
  equation PS-3 are split with respect to the reflection
  symmetry into two groups: the {\it  symmetric} and the {\it
  antisymmetric}. In the present paper we consider only the second group of
  solutions.

\begin{thrm}
  When $\lambda\neq 1,3$ the antisymmetric eigenfunctions $u(x)$ of
  the PS-3
  integral equation with parameter $R_3(x)$ are in one to one correspondence with the pants
  ${\cal P}_s(\lambda, h_1,h_2|m_1,m_2)$, $s=1,2,3,12,13$, which are conformally
  equivalent to the pair of pants ${\cal P}(R_3)$ with colored boundary components.

  Let $p(y)$ be the conformal map from ${\cal P}(R_3)$  to ${\cal P}_s(\lambda, h_1,h_2|m_1,m_2)$, then
  up to proportionality

\begin{equation}
u(x)=
\sqrt{\frac{(y-y_1)(y-y_2)}{p'(y^+)p'(y^-)}}
\frac{p(y^+)-p(y^-)}{w(y)},
\label{UviaP}
\end{equation}
where $x\in[-1,1]$;  $y:=R_3(x)$, $y^\pm:=y\pm i0$. For $s=1$,
$y_1=\overline{y_2}$ is the critical point of the function $p(y)$;
for $s=2,3$ the real $y_1$ and $y_2$ are critical points of the function
$p(y)$.
\end{thrm}

The p~r~o~o~f of this theorem will be given in the remaining two sections of the article.

\subsection{Corollaries}
The representation (\ref{UviaP}) cannot be called explicit in the usual sense,
since it comprises a transcendent function $p(y)$. We show that nevertheless
the representation is useful as it allows us to understand the following
properties of the solutions.

{\bf 1.} The "antisymmetric" part of the spectrum is always a subset of $[1,2)\cup\{3\}$.

{\bf 2.} Every $\lambda\in (1,2)$ is the eigenvalue for infinitely many equations PS-3.

P~r~o~o~f.~ Any of the constructed pants may be transformed to the
standard  form: a sphere with three real slots. Normalizing the red
slot to be $[-1,1]$, the end points of the two other slots will give the
branch points $a_1$, \dots, $a_4$. We know already how to
reconstruct the branched covering $R_3(x)$ from its branch points.

{\bf 3.} Eigenfunction $u(x)$ related to the pants ${\cal P}_s(\dots|m_1,m_2)$
has exactly $m_1+m_2+2$ zeroes on the segment $[-1,1]$ when $s=2,3$ and one more zero when $s=1$.

  P~r~o~o~f.~ According to the formula
  (\ref{UviaP}), the number of zeroes of eigenfunction  $u(x)$ is equal to
  the number of points $y\in[-1,1]$ where $p(y^+)=p(y^-)$. This
  number in turn is equal to the number of solutions of the
  inclusion

\begin{equation}
S(y):=
Arg [p(y^-)-\mu^{-1}] -
Arg [p(y^+)-\mu^{-1}]
\quad\in 2\pi\mbox{\sets Z},
\qquad y\in[-1,1].
\label{Inclus}
\end{equation}

  Let the point $p(y)$ go $m$ times around the circle $C$  when
  its argument $y$ travels along the two sides of $[-1,1]$. The integer $m$
  is naturally related to the integer parameters of the pants ${\cal
  P}_s(\dots)$. The function $S(y)$ strictly increases from $0$
  to $2\pi m$ on the segment $[-1,1]$, therefore the inclusion
  (\ref{Inclus}) has exactly $m+1$ solutions on the mentioned
  segment.

  {\bf 4.} The mechanism for generating the discrete spectrum of the integral
  equation is explained. Sewing an annulus to the pants ${\cal
  P}_s(\lambda,h_1,h_2|\dots)$ changes the conformal structure of
\index{pants!conformal structure of}
  the latter. To return to the conformal structure specified by ${\cal P}(R_3)$ we have to
  change the real parameters of the pants, one of them being the
  spectral parameter $\lambda$.

If we knew how to evaluate the conformal moduli of the pair of pants
${\cal P}_s(\lambda, h_1,h_2|m_1,m_2)$ as functions of its real
parameters, the solution of the integral equation would be reduced
to a system of three transcendental equations for the three numbers
$\lambda, h_1,h_2$. This solution will  depend on the integer
parameters $s$, $m_1$, $m_2$.


\section{Geometry of Integral Equation}
PS integral equations possess a rich geometrical structure which we disclose in this section.
The chain of equivalent transformations of PS-3 equation described
here in a somewhat sketchy fashion
is given in \cite{Bog2, Bog3} with more details.

\subsection{A nonlocal Functional Equation}\label{NonFunEq}
\index{equation!functional} Let us decompose the kernel of the
second integral in (\ref{PSE}) into a sum of elementary fractions:
\begin{equation}
\label{kernel}
\frac{ R_3'(t)}{R_3(t)-R_3(x)}=
\frac{d}{dt}\log(R_3(t)-R_3(x))=
\sum_{k=1}^3\frac{1}{t-x_k(x)}-\frac{Q'}{Q}(t),
\end{equation}
where $Q(t)$ is the denominator in an irreducible representation of
$R(t)$ as the ratio of two polynomials; $x_1(x)=x,~x_2(x),~x_3(x)$
-- are all solutions (including multiple and infinite) of the
algebraic equation $R_3(x_s)=R_3(x)$. This expansion suggests to
rewrite the original equation (\ref{PSE}) as a certain relationship
for the Cauchy-type integral
\begin{equation}
\Phi(x):=\int\limits_I\frac{u(t)}{t-x}dt+const^*,
\quad x\in\hat{\mbox{\sets C}}\setminus[-1,1].
\label{Cauchy}
\end{equation}
The constant term $const^*$ in (\ref{Cauchy}) is introduced to compensate for the constant terms arising after
substitution of expression (\ref{Cauchy}) to the equation (\ref{PSE}).

For a known $\Phi(x)$, the eigenfunction $u(t)$ may be recovered by the So\-khot\-skii-Ple\-melj
formula:
\begin{equation}
\label{SP}
u(t)=(2\pi i)^{-1}
\left[
\Phi(t+i0)-\Phi(t-i0)
\right],
\quad t\in I.
\end{equation}

Function $\Phi(x)$ generated by an eigenfunction of PS integral equation satisfies a nonlocal functional equation:

\begin{lmm}~\cite{Bog2}~~~
For $\lambda\neq 1,3$ the transformations (\ref{Cauchy}) and
(\ref{SP}) imply a
1-1 correspondence between the H\"older eigenfunctions $u(t)$ of the PS-3
integral equation and the nontrivial solutions $\Phi(x)$ of the
functional equation which are holomorphic on a sphere with the slot $[-1,1]$
\begin{equation}
\label{FE}
\Phi(x+i0)+\Phi(x-i0)=\delta
\biggl(
\Phi (x_2(x))+\Phi(x_3(x))
\biggl)~,
\quad x\in I,
\end{equation}
\begin{equation}
\label{delta}
\delta =2/(\lambda -1),
\end{equation}
with H\"older boundary values $\Phi(x\pm i0)$.
\end{lmm}

\subsection{The Riemann Monodromy Problem}
\label{RiMoPr}\index{Riemann!monodromy problem}

  The lifting $R_3^{-1}({\cal P}(R_3))$ of the pants associated
  to the integral equation consists of three components ${\cal
  O}_s$,~~ $s=1,2,3$. We number them in the following way (see Fig.
  \ref{CoverTopology}): the segment $[-1,1]$ lies on the boundary
  of ${\cal O}_1$; the segment $[c_4,c_3]$ is on the boundary of
  ${\cal O}_2$ and the boundary of ${\cal O}_3$ comprises the
  segment $[c_2,c_1]$.

\subsubsection{}
\begin{figure}[ht!]
\centerline{\psfig{figure=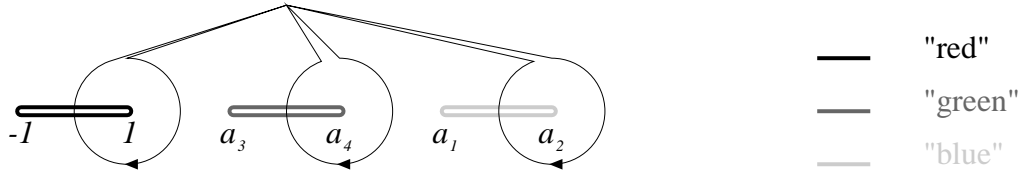}} \caption[] {\normalsize
Three loops on a sphere with six punctures $\pm1$, $a_1,\dots$,
$a_4$} \label{ThreeLoops}
\end{figure}

  Let the function $\Phi(x)$ be related to the solution $u(x)$ of the integral equation (\ref{PSE}) as in (\ref{Cauchy}).
  We consider a 3-vector defined in the pair of pants:

\begin{equation}
\label{Wofy}
W(y):=(W_1,W_2,W_3)^t=(\Phi(x_1),\Phi(x_2),\Phi(x_3))^t,
\qquad y\in{\cal P}(R_3),
\end{equation}
  where  $x_s$ is the unique
  solution of the equation $R_3(x_s)=y$ in ${\cal O}_s$. Vector
  $W(y)$ is holomorphic and bounded in the pants ${\cal P}(R_3)$ as all three points $x_s$, $s=1,2,3$,
  remain in the holomorphicity domain of the function $\Phi(x)$. We claim
  that {\it the boundary values of the vector $W(y)$ are related
  via constant matrices}:
\begin{equation}
\label{RMPr}
  W(y+i0)={\bf D_*}W(y-i0),
  \qquad ~~{\rm when}~~ y\in\{slot_*\}.
\end{equation}
The matrix ${\bf D_*}$ assigned
 to the  "green" $[a_3,a_4]$, "blue" $[a_1,a_2]$, and "red" $ [-1,1]$  slot respectively is
\begin{equation}
\label{permutations}
{\bf D_2}:=
\begin{array}{||lll||}
0&0&1\\0&1&0\\1&0&0\\
\end{array}~;
\qquad
{\bf D_3}:=
\begin{array}{||lll||}
0&1&0\\1&0&0\\0&0&1\\
\end{array}~;
\qquad
{\bf D}:=
\begin{array}{||ccc||}
-1       & \delta & \delta \\
 0       & 1      & 0  \\
 0       & 0      & 1  \\
\end{array}
~~.
\end{equation}
  This in particular means that our vector (\ref{Wofy}) is a
  solution of a certain Riemann monodromy problem. The monodromy
  of vector $W(y)$ along the loop crossing only "red", "green" or
  "blue" slot is given by the matrix $\bf D$, $\bf D_2$ or $\bf D_3$
  correspondingly -- see Fig. \ref{ThreeLoops}.
  \index{monodromy!generators}

  Indeed, let $y^+:=y+i0$ and $y^-:=y-i0$ be two points on opposite
  sides of the "blue" slot $[a_1,a_2]$. Their inverse images $x_3^+=x_3^-$, $x_1^\pm=x_2^\mp$
  lie outside the cut $[-1,1]$. Hence $W(y^+)=$
  ${\bf D_3}W(y^-)$. For two points $y^\pm$ lying on opposite
  sides of the "green" slot $[a_3,a_4]$, their inverse images satisfy
  relations $x_2^+=x_2^-$, $x_1^\pm=x_3^\mp$, which means
  $W(y^+)=$ ${\bf D_2}W(y^-)$. Finally, let $y^\pm$ lie on both  sides
  of the "red" slot $[-1,1]$. Now two points $x_2^+=x_2^-$ and $x_3^+=x_3^-$ lie
  in the holomorphicity domain of $\Phi(x)$ while $x_1^+$ and $x_1^-$
  appear on the opposite sides of the cut $[-1,1]$. According to
  the functional equation (\ref{FE}), \index{slot}
$$
\Phi(x_1^+)=
-\Phi(x_1^-)+\delta(\Phi (x_2^\pm)+\Phi (x_3^\pm)),
$$
therefore $W(y^+)=$ ${\bf D}W(y^-)$ holds on the slot $[-1,1]$.

\subsubsection{}
  Conversely, let $W(y)$ be the bounded solution of the Riemann monodromy
  problem (\ref{RMPr}). We define a piecewise holomorphic
  function  on the Riemann sphere:
  \begin{equation}
  \Phi(x):=W_s(R_3(x)),
  \qquad {\rm when}~~x\in {\cal O}_s,
  \qquad s=1,2,3.
  \label{PhifromW}
  \end{equation}
  From the boundary relations for the vector $W(y)$ it
  immediately follows that the function $\Phi(x)$ has no jumps on
  the lifted cuts $[a_1,a_2]$, $[a_3,a_4]$, $[-1,1]$
  except for the cut $[-1,1]$ from the upper sphere. Say, if the two points
  $y^\pm$ lie on opposite sides of the cut $[a_1,a_2]$, then
  $W_3(y^+)=W_3(y^-)$ and  $W_1(y^\pm)=W_2(y^\mp)$ which means
  that the function $\Phi(x)$ has no jump on the components of
  $R_3^{-1}[a_1,a_2]$. From the boundary relation on the cut $[-
  1,1]$ it follows that $\Phi(x)$ is the solution for the functional
  equation (\ref{FE}). Therefore it gives a solution of Poincar\'{e}--Steklov
  integral equation with parameter $R_3(x)$. Combining
  formulae (\ref{SP}) with (\ref{PhifromW}) we get the
  reconstruction rule
\begin{equation}
\label{uofx}
u(x)=(2\pi i)^{-1}\biggl(W_1(R_3(x)+i0)-W_1(R_3(x)-i0)\biggr),
\qquad x\in [-1,1].
\end{equation}

We have just proved the following
\begin{thrm} \cite{Bog1}
~~If $\lambda\neq1,3$ then the two formulas (\ref{Wofy}) and (\ref{uofx})
imply the one-to-one  correspondence between the solutions $u(x)$ of
the integral equation (\ref{PSE}) and the bounded solutions $W(y)$ of the
Riemann monodromy problem (\ref{RMPr}) in the punctured sphere $\hat{\mbox{\sets C}}\setminus\{
-1,1,a_1,a_2,a_3,a_4\}$.
\end{thrm}

\subsubsection{Monodromy Invariant}
\index{monodromy!invariant}
The following statement is proved by direct computation.
\begin{lmm}
All matrixes (\ref{permutations}) (i) are involutive (i.e. ${\bf
D}^2={\bf D_2}^2={\bf D_3}^2={\bf 1}$) and (ii) conserve the
quadratic form
\begin{equation}
\label{J}
J(W):=\sum\limits_{k=1}^3 W_k^2-
\delta\sum\limits_{j<s}^3 W_jW_s.
\end{equation}
\end{lmm}

  The form $J(W)$ is not degenerate unless $-2\ne\delta\ne
  1$, or equivalently $0\ne\lambda\ne 3$. Since the solution
  $W(y)$ of our monodromy problem is bounded near the cuts, the
  value of the form $J(W)$ is independent of the variable $y$.
  Therefore the solution takes values either in the smooth quadric
  $\{W\in\mbox{\sets C}^3:~~J(W)=J_0\neq0\}$, or the cone
  $\{W\in\mbox{\sets C}^3:~~J(W)=0\}$.

\subsection{Geometry of the Quadric Surface}
\begin{figure}
\centerline{\psfig{figure=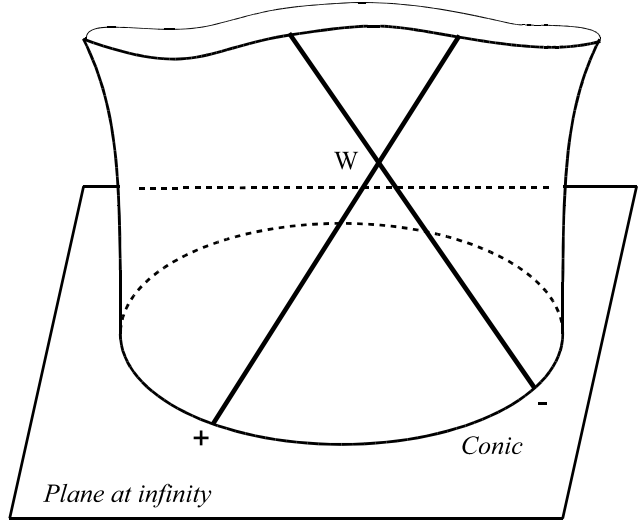}} \caption[Quadric] {\small
Global coordinates $p^+$ and $p^-$ on quadric } \label{Quadric}
\end{figure}

\index{quadric}\index{quadric!line generators}
The nondegenerate projective quadric $\{J(W)=J_0\}$ contains two families of
line elements\footnote{
This property of quadric is sometimes used in architecture. The line generators of
the hyperboloid serve as construction elements, e.g.,  for the
Shukhov tower in Moscow.}
which for convenience are denoted by the signs $'+'$ and $'-'$.
Two different lines from the same family are disjoint while two lines from
different families must intersect. The intersection of those lines with
the 'infinitely distant' secant plane gives points on the conic
\index{conic}
\begin{equation}
\label{conica}
{\cal C}:=\{(W_1:W_2:W_3)^t\in\mbox{\sets CP}^2:\quad J(W)=0\}
\end{equation}
which by means of the stereographic projection $p$ may be identified
with the Riemann sphere. Therefore we have introduced two global
coordinates $p^\pm(W)$ on the quadric, the 'infinite part' of which (=
conic $\cal C$) corresponds to coinciding coordinates: $p^+=p^-$
(see fig. \ref{Quadric}).

The natural action of them pseudo-orthogonal group $O_3(J)$ in $\mbox{\sets C}^3$
conserves the quadric, the conic at infinity $\cal C$, and the families of line elements
possibly interchanging their labels $'\pm'$.  The induced action of the group
$O_3(J)$ on the stereographic coordinates $p^\pm$ is a linear
fractional with a possible change
of the superscript $'\pm'$.

\subsubsection{Stereographic Coordinates}
\index{stereographic coordinates}
To obtain explicit expressions for the coordinate change
$W\leftrightarrow p^\pm$ on the quadric
we bring the quadratic form $J(W)$ to the simpler form
$J_\bullet(V):=$ $V_1V_3-V_2^2$ by means of the linear coordinate change
\begin{equation}
W={\bf K}V
\label{V}
\end{equation}
where
\begin{equation}
\label{K}
{\bf K}:=
(3\delta+6)^{-1/2}
\begin{array}{||lll||}
1&1&1\\
1&\varepsilon^2&\varepsilon\\
1&\varepsilon&\varepsilon^2
\end{array}
\cdot
\begin{array}{||lll||}
0&\mu^{-1}&0\\
0&0&1\\
1&0&0
\end{array}
~,
\end{equation}
$$
\varepsilon:=\exp(2\pi i/3),
\qquad
\mu:=\sqrt{\frac{\delta-1}{\delta+2}}=
\sqrt{\frac{3-\lambda}{2\lambda}}.
$$
Translating the first paragraph of the current section
into formulae we get
\begin{equation}
\label{V2Ppm}
p^\pm(W):=\frac{V_2\pm i\sqrt{J_0}}{V_1}
=\frac{V_3}{V_2\mp i\sqrt{J_0}};
\end{equation}
and inverting this dependence,
\begin{equation}
\label{Ppm2W}
W(p^+,p^-)=\frac{2i\sqrt{J_0}}{p^+-p^-}{\bf K}
\left(
\begin{array}{c}
1\\
(p^++p^-)/2\\
p^+p^-
\end{array}
\right)~.
\end{equation}
The point $W(p^+,p^-)$ with coordinate $p^+$ (resp. $p^-$) being
fixed moves on the straight
line with the directing vector ${\bf K}(1:p^+:(p^+)^2)$ (resp. ${\bf K}(1:p^-:(p^-)^2)$)
belonging to the conic (\ref{conica}).

\subsubsection{Action of the Pseudo-orthogonal Group}
\index{spinor representation}
\begin{lmm}
There exists a (spinor) representation $\chi:\quad O_3(J)\to PSL_2(\mbox{\sets C})$ such that:\\
1) The restriction of $\chi(\cdot)$ to $SO_3(J)$ is an isomorphism to $PSL_2(\mbox{\sets C})$.\\
2) For coordinates $p^\pm$ on the quadric the following transformation rule holds:\\
\begin{equation}
\label{Ptrans}
\begin{array}{ll}
p^\pm({\bf T}W)=\chi({\bf T})p^\pm(W),&
\qquad{\bf T}\in SO_3(J),\\
p^\pm({\bf T}W)=\chi({\bf T})p^\mp(W),&
\qquad{\bf T}\not\in SO_3(J).\\
\end{array}
\end{equation}
3) The linear-fractional mapping $\chi p:= (ap+b)/(cp+d)$ is the image of the matrix:
\begin{equation}
\label{chiT}
{\bf T}:=\frac1{ad-bc}~{\bf K}~
\begin{array}{||ccc||}
d^2&2cd&c^2\\
bd&ad+bc &ac\\
b^2&2ab&a^2
\end{array}
~{\bf K}^{-1}
\in SO_3(J).
\end{equation}
4) The generators of the monodromy group are mapped to the following
elements of $PSL_2$:
\begin{equation}
\label{chiD}
\begin{array}{ll}
\chi({\bf D}_s)p=\varepsilon^{1-s}/p,&\qquad s=1,2,3;\\
\chi({\bf D})p=
\displaystyle{\frac{\mu p-1}{p-\mu}}.&
\end{array}
\end{equation}
\end{lmm}
P~r~o~o~f: We define the action of the matrix ${\bf A}\in
SL_2(\mbox{\sets C})$ on the vector $V\in\mbox{\sets C}^3$ by the
formula:

\begin{equation}
\label{SL2action}
{\bf A}:=\begin{array}{||cc||}a&b\\c&d\end{array}:
\qquad
\begin{array}{||cc||}
V_3&V_2\\
V_2&V_1
\end{array}
\quad
\longrightarrow
\quad
{\bf A}~
\begin{array}{||cc||}
V_3&V_2\\
V_2&V_1
\end{array}
~{\bf A}^t.
\end{equation}
  It is easy to check that (\ref{SL2action}) gives the faithful
  representation of a connected 3-dimensional group
  $PSL_2(\mbox{\sets C}):=SL_2(\mbox{\sets C})/\{\pm {\bf 1}\}$ into
  $SO_3(J_\bullet)$ and therefore, an isomorphism. Let us  denote by
  $\chi_\bullet$ the inverse isomorphism  $SO_3(J_\bullet)\to$
  $PSL_2(\mbox{\sets C})$ and let $\chi(\pm {\bf
  T}):=\chi_\bullet({\bf K}^{-1}{\bf TK})$ for ${\bf T}\in
  SO_3(J)$. The obtained homomorphism $\chi:~~$ $O_3(J)\to$
  $PSL_2(\mbox{\sets C})$ will satisfy statement 1) of the lemma.

  To prove 2) we replace components of the vector $V$ in the right-hand
  side of (\ref{SL2action}) with their representation in terms of
  the stereographic coordinates $p^\pm$:

$$
\frac{i\sqrt{J_0}}{p^+-p^-}
{\bf A}~
\left[
(p^+,1)^t\cdot(p^-,1)+
(p^-,1)^t\cdot(p^+,1)
\right]
~{\bf A}^t=
$$
$$
i\sqrt{J_0}
\frac{(cp^++d)(cp^-+d)}{p^+-p^-}
\left[
(\chi p^+,1)^t\cdot(\chi p^-,1)+
(\chi p^-,1)^t\cdot(\chi p^+,1)
\right]
=
$$
$$
\frac{i\sqrt{J_0}}{\chi p^+-\chi p^-}
\left[
(\chi p^+,1)^t\cdot(\chi p^-,1)+
(\chi p^-,1)^t\cdot(\chi p^+,1)
\right]
=
$$
$$
\begin{array}{||cc||}
V_3(\chi p^+,\chi p^-)&V_2(\chi p^+,\chi p^-)\\
V_2(\chi p^+,\chi p^-)&V_1(\chi p^+,\chi p^-)
\end{array}
~~,
$$
  where we set $\chi p:=(ap+b)/(cp+d)$. Now (\ref{Ptrans})
  follows immediately for ${\bf T}\in SO_3(J)$. It remains to
  check the transformation rule for any matrix ${\bf T}$ from
  the other component of the group $O_3(J)$, say ${\bf T}=-{\bf 1}$.

  Writing the action (\ref{SL2action}) component-wise we arrive
  at conclusion 3) of the lemma.

  And finally, expressions 4) for the generators of monodromy
  group may be obtained either from analyzing formula
  (\ref{chiT}) or from finding the eigenvectors of the matrices ${\bf
  D}_s,{\bf D}$ which correspond to the fixed points of linear-
  fractional transformations. ~~~\bl

For convenience we collect all the introduced objects related to the
boundary components of the pair of pants ${\cal P}(R_3)$ in
table \ref{Tab1}
\begin{table}
\label{Tab1} \caption{Slots, their associated colors, matrices and
linear-fractional maps}
\centering
\begin{tabular}{c||c|c|c}
slot& $[-1,1]$&$[a_1,a_2]$&$[a_3,a_4]$\\
\hline
color& "red" & "blue" & "green"  \\
\hline
matrix $\bf D_*$&
$\displaystyle{
{\bf D}:=
\begin{array}{||ccc||}
-1       & \delta & \delta \\
 0       & 1      & 0  \\
 0       & 0      & 1  \\
\end{array}}
$&
$\displaystyle{
{\bf D}_3:=
\begin{array}{||lll||}
0&1&0\\1&0&0\\0&0&1\\
\end{array}}
$&
$\displaystyle{
{\bf D}_2:=
\begin{array}{||lll||}
0&0&1\\0&1&0\\1&0&0\\
\end{array}}
$\\
\hline
$\chi({\bf D_*})p$&
$\displaystyle{\frac{\mu p-1}{p-\mu}}$&
$\varepsilon/p$&$\varepsilon^2/p$
\end{tabular}
\end{table}

\subsection{Entangled Projective Structures}\label{CouProjStr}
\index{projective strucure!branched}\index{projective structure!monodromy of}

  {\bf Definition.} A branched complex projective structure \cite{G,Man,Tu,GKM} on a Riemann
  surface $\cal M$ is a meromorphic function $p(t)$ on the
  universal covering $\widetilde{\cal M}$ which transforms fractionally linear
  under the cover transformations of
  $\widetilde{\cal M}$. The appropriate representation
  $\chi:\quad$ $\pi_1({\cal M})\to $ $PSL_2(\mbox{\sets C})$ is called
  the {\it  monodromy} of the projective structure. The set of
  all critical points of $p(t)$ with their multiplicities survives
  under the cover transformations of $\widetilde{\cal M}$.
  The projection of this set to the Riemann surface ${\cal M}$ is known as
  the  {\it  branching divisor} ${\sf D}(p)$ of projective structure
  and the branching number of the structure $p(t)$ is $\deg{\sf D}(p)$.

\index{projective strucure!unbranched} {\bf Examples.} The
unbranched projective structures arise in Fuchsian and Schottky
uniformizations of the Riemann surface. Any meromorphic function on
a Riemann surface is a branched projective structure with trivial
monodromy. \index{uniformization}

\subsubsection{Projective Structures Generated by Eigenfunction}
Every bounded solution $W(y)$
of the Riemann monodromy problem (\ref{RMPr}) generates
two nowhere coinciding  meromorphic functions $p^\pm(y)$ in the sphere
with three slots.  Those functions are stereographic coordinates (\ref{V2Ppm}) for
the vector $W(y)$. The boundary values of functions $p^+(y)$ and $p^-(y)$ on every slot
are related by linear-fractional transformations:
\begin{equation}
\label{BVp}
p^\pm(y+i0)=\chi({\bf D}_*)p^\mp(y-i0),
\qquad y\in \{slot_*\}
\end{equation}
where the matrix ${\bf D}_*=$ ${\bf D}$, ${\bf D}_2$, ${\bf D}_3$
stand
for the 'red', 'green' and 'blue' slots respectively.

\index{universal covering}
Relations (\ref{BVp}) allow us to analytically continue both functions $p^+(y)$ and $p^-(y)$
through any slot to the second sheet of the genus 2 Riemann surface
\begin{equation}
M:=\{w^2=(y^2-1)\prod\limits_{s=1}^4(y-a_s)\},
\label{M}
\end{equation}
and further to its universal covering $\tilde{M}$. Thus obtained
functions $p^\pm(t)$, $t\in\tilde{M}$, will be locally single valued
on the Riemann surface since all matrices ${\bf D}_*$ are
involutive. However varying the argument $t$ along the handle
of the surface $M$ may result in a linear-fractional
transformation of the value $p^\pm(t)$. Say, the continuations of
$p^+(y)$ from the pants through the red and green slots will give
two different functions on the second sheet related by the
linear-fractional mapping $\chi({\bf DD}_2)$.

\subsubsection{Branching of the Structures $p^\pm$}
The way we have carried out the continuation of functions $p^\pm(y)$ suggests
that the branching divisors of the arising projective structures are related via
the hyperelliptic involution $H(y,w):=(y,-w)$ of the surface $M$:
\begin{equation}
\label{BranchSymm}
{\sf D}(p^+)=H{\sf D}(p^-).
\end{equation}
The condition $p^+\neq p^-$ allows to determine the branching numbers of the structures
which is done in the next theorem.

\begin{thrm} \cite{Bog3} \label{J0neq0}
  When $\lambda\not\in\{0,1,3\}$ the solutions $u(x)$ of the  integral equation PS-3
  that have invariant $J_0\neq0$ are in
  one-to-one correspondence with the couples of not identically
  equal functions meromorphic in the pants ${\cal P}(R_3)$
  $p^\pm(y)$ with boundary values satisfying (\ref{BVp}) and two
  critical points in common. The correspondence $u(x)\to$
  $p^\pm(y)$ is established by the sequence of formulae
  (\ref{Cauchy}), (\ref{Wofy}) and (\ref{V2Ppm}); the inverse
  dependence is given by the formula
\begin{equation}
\label{UfromP}
2\pi u(x)=\sqrt{\frac{(\delta+2)J_0}3}
\frac{p^+(y)p^-(y) -\mu(p^+(y)+p^-(y))+1}
{p^+(y)-p^-(y)},
\end{equation}
where $x\in[-1,1]$ and $y:=R_3(x)+i0$.
\end{thrm}
{\bf Remark}:  The number of critical points of the structures in
the pants is counted with their {\it  weight and multiplicity}
(see remark 3 on page \pageref{Rmk3}): ~~1) the branching number
of $p^\pm(y)$ at the
  branch point  $a\in \{\pm1, a_1,\dots,a_4\}$ of $M$ is computed
  with respect to the local parameter $z=\sqrt{y-a}$, ~~2) every
  branch point of the projective structure on the boundary of the
  pants should be considered as a half-point.

{\bf P~r~o~o~f}: {\bf 1.~~} Let $u(x)$ be an eigenfunction of the
integral equation PS-3, then the stereographic coordinates
$p^\pm(y)$ of the solution of the associated Riemann monodromy
problem inherit the boundary relationship (\ref{BVp}). What
remains is to find the branching numbers of the entangled
structures $p^\pm(y)$. To this end we consider the {\it Kleinian
quadratic differential} on the slit sphere \index{quadratic
differential!Kleinian}
\begin{equation}
\label{Klein}
\Omega(y)=
\frac{dp^+(y)dp^-(y)}
{(p^+(y)-p^-(y))^2},
\qquad y\in\hat{\mbox{\sets C}}.
\end{equation}
This expression is the infinitesimal form of the cross ratio, hence
it remains unchanged after the same linear-fractional
transformations of the functions $p^+$ and $p^-$. Therefore,
(\ref{Klein}) is a well defined quadratic differential on the entire
sphere. Lifting $\Omega(y)$ to the surface $M$ we get a holomorphic
differential. Indeed, $p^+\neq p^-$ everywhere and applying suitable
linear-fractional transformation we assume that
$p^+=1+z^{m_+}+\{terms~of~higher~order\}$ and $p^-=cz^{m_-}+...$ in
terms of local parameter $z$ of the surface, $m_\pm\ge1$, $c\neq0$.
Then $\Omega=cm_+m_-z^{m_++m_--2}+\{terms~of~higher~order\}$.
Therefore
$$
{\sf D}(p^+)+{\sf D}(p^-)=(\Omega).
$$
Any holomorphic quadratic differential on a genus 2 surface has 4 zeroes and taking
into account the symmetry (\ref{BranchSymm}) of the branching divisors,
we see that each of the structures $p^\pm$ has the branching number two on the curve $M$.
It remains to note that the pair of pants ${\cal P}(R_3)$ are exactly "one half" of $M$.

{\bf 2.~~}
Conversely, let $p^+(y)$ and $p^-(y)$ be two not identically equal meromorphic functions on the slit sphere,
with boundary conditions (\ref{BVp}) and total branching number two in the pants (see remark above). We can prove that $p^+\neq p^-$ everywhere.
Indeed, for the meromorphic quadratic differential (\ref{Klein}) on the Riemann surface $M$
we establish (using a local coordinate on the surface) the inequality
\begin{equation}
\label{DivIneq}
{\sf D}(p^+)+{\sf D}(p^-)\ge(\Omega)
\end{equation}
where the deviation from equality means that there is a point where $p^+=p^-$. But
the degree of the divisor on the left of (\ref{DivIneq}) is four and the same number is $\deg(\Omega)=4g-4$.
Therefore this pair of functions $p^\pm$ will give us the holomorphic vector $W(p^+(y),p^-(y))$
in the pants which solves our Riemann monodromy problem. We already know how to convert the latter vector to the
eigenfunction of the integral equation PS-3. ~~~\bl

\subsubsection{Remark about the Non-smooth Quadric}\label{ConeCase}

  It is shown in \cite{Bog3} how to
  incorporate the exceptional case $J_0=0$ into the above scheme.
  In the latter case the functions $p^\pm(y)$ coincide, however the
  boundary relations (\ref{BVp}) survive. The total
  branching number of the function $p^+=p^-$ in the pair of pants is
  either zero or one. The solutions to the PS-3
  integral equation and the associated Riemann monodromy problem may be recovered up
  to proportionality from the unified formulae (true whatever
  $J_0$)
\begin{equation}
\label{UniU}
u(x)=\sqrt{\frac{\Omega(y)}{dp^+(y)dp^-(y)}}
(p^+(y)p^-(y) -\mu(p^+(y)+p^-(y))+1),
\end{equation}
\begin{equation}
\label{UniW}
W(y)=\sqrt{\frac{\Omega(y)}{dp^+(y)dp^-(y)}}
{\bf K}(1, (p^+(y)+p^-(y))/2, p^+(y)p^-(y))^t,
\end{equation}
  where $\displaystyle{\Omega(y)=(y-y_1)(y-y_2)\frac{(dy)^2}{w^2(y)}}$
  is the holomorphic quadratic differential on
  the Riemann surface $M$ with zeroes at the branching points of
  the possibly coinciding structures $p^+$ and $p^-$ [or with
  two arbitrary double zeroes when the structure $p^+=p^-$ is unbranched,
  further analysis however shows that the required unbranched structures do not exist].

\subsection{Types of the Mirror Symmetry of the  Solution}
The eigenvalues of the integral equation are the critical values
of the {\it  positive} functional (\ref{Rayleigh}) -- the
generalized Rayleigh ratio. So we may consider only {\it  real}
eigenfunctions $u(x)$ without loss of generality. Real solutions
of the PS-3 equation give rise to exactly two types of {\it
mirror symmetry} for the
  entangled structures:
$$
\begin{array}{ll}
Symmetric& p^\pm(\bar{y})=1/\overline{p^\pm}(y)\\
Antisymmetric& p^\pm(\bar{y})=1/\overline{p^\mp}(y) \\
\end{array},
\qquad y\in {\cal P}(R_3),
$$
depending on the sign of the real number $(\delta+2)J_0$. In what
follows we restrict ourselves to the case of {\it  antisymmetric
eigenfunctions}. In this case:
$$
p^+(y\pm i0)=
1/\overline{p^-(y\mp i0)}=
1/\overline{\chi({\bf D_*})p^+(y\pm i0)},
\qquad y\in slot_*,
$$
and hence we know where the boundary components of the pair of pants
${\cal P}(R_3)$ are mapped to. In particular,
\begin{equation}
\begin{array}{ll}
{\rm "green"~boundary} & \to\varepsilon\hat{\mbox{\sets R}}\\
{\rm "blue"~boundary}&   \to \varepsilon^2\hat{\mbox{\sets R}}\\
{\rm "red"~boundary}&    \to \left\{
\begin{array}{ll}
C \quad - see~ (\ref{circle}) & ~when~ 1<\lambda<3 \\
\emptyset & ~when~~ \lambda<1 ~or~ 3<\lambda.
\end{array}
\right.
\end{array}
\label{MapBoundary}
\end{equation}
We see that the above geometrical analysis of the integral equation gives
the universal limits for (the antisymmetric part of) the spectrum.

  The branching divisor of the projective structure $p^+$ has
  the mirror symmetry: ${\sf D}(p^+)=\bar{H} {\sf D}(p^+)$ where
  $\bar{H}(y,w):=(\bar{y},-\bar{w})$ is the anticonformal
  involution of the surface $M$ leaving boundary components (ovals) of
  pair of the pants ${\cal P}(R_3)$ intact. Therefore exactly three situations may occur: $p^+(y)$ has
  one simple critical point strictly inside the pants, or
  there are two simple critical points on the boundary of
  pants or  there is one double critical point
  of $p^+(y)$ on the boundary of the pants.

\section{Combinatorics of Integral Equation}
\label{ProblemPants}
For the antisymmetric eigenfunctions we arrive at the
essentially combinatorial \\
{\bf Problem} (about putting pants on a sphere)~~
  {\it Find a meromorphic function $p:=p^+$ defined in the pair of pants ${\cal P}(R_3)$
  mapping boundary ovals to the given circles (\ref{MapBoundary}) and having
  exactly one critical point (counted with weight and multiplicity) in the pants.}

  The three above-mentioned types of the branching divisor ${\sf
  D}(p)$ will be treated separately in the sections \ref{Inside}, \ref{OnBorder}.
  When the branch point of the structure $p$ is
  strictly inside the pants we show that the solution of the problem
  takes  the form of the \"Uber\-la\-ge\-rungs\-fl\"a\-che ${\cal
  P}_1(\dots)$ with certain real and integer parameters. The case
  of two simple branch points belonging to the boundary gives us
  the pants ${\cal P}_s(\dots)$, $s=2,3$ and the unstable
  intermediate case with double branch point on the boundary
  brings us to the pants ${\cal P}_j(\dots)$, $j=12,~13$
  described in (\ref{InterMediate}). \index{\"Uberlagerungsfl\"ache}

  Let $p(y)$ be a holomorphic map from a Riemann surface $\cal M$
  with a boundary  to the sphere and the selected boundary
  component $(\partial {\cal M})_*$ be mapped to a circle. The
  reflection principle allows us to holomorphically continue
  $p(y)$ through this selected component to the double of $\cal
  M$. Therefore we can talk of the critical points of $p(y)$ on
  $(\partial {\cal M})_*$. When the argument $y$ passes through a
  simple critical point, the value $p(y)$ reverses the direction
  of its movement on the circle. So there should be an even number
  of critical points (counted with multiplicities) on the
  selected boundary component.

\subsection{The Branchpoint is inside a Pair of Pants}\label{Inside}
\subsubsection{Construction 1}\label{Construct1}
Using otherwise a composition with a suitable linear-fractional map, we suppose that the circle
$p((\partial {\cal M})_*)$ is the boundary of the unitary disc
\begin{equation}
\label{Udisc}
\mbox{\sets U}:=\{p\in\mbox{\sets C}:\quad |p|\le 1\},
\end{equation}
and that a small annular vicinity of the selected boundary component is
mapped to the exterior of the unit disc. We define the mapping of a
disjoint union ${\cal M}\cup\mbox{\sets U}$ to a sphere
\begin{equation}
\tilde{p}(y):=\left\{
\begin{array}{ll}
p(y),& y\in {\cal M},\\
L(y^d),& y\in \mbox{\sets U},
\end{array}
\right.
\end{equation}
where the integer $d>0$ is the degree of the mapping
$p:$  $(\partial {\cal M})_*\to$ $\partial${\sets U},
and where $L(y)$ is an (at the moment arbitrary) linear fractional mapping
keeping the unitary disc (\ref{Udisc}) unchanged. The choice of $L(\cdot)$
will be fixed later to simplify the arising combinatorial analysis.
\index{winding number}

  Now we fill in the hole in ${\cal M}$ by the unit disc,
  identifying the points of $(\partial {\cal M})_*$ and the
  points of $\partial${\sets U} with the same value of
  $\tilde{p}$ (there are $d$ ways to do so). The holomorphic
  mapping $\tilde{p}(y)$ of the new Riemann surface ${\cal
  M}\cup${\sets U} to the sphere will have exactly one additional
  critical point of multiplicity $d-1$ at the center of the glued
  disc.

\subsubsection{Branched Covering of a Sphere}
  We return to the function $p(y)$ being the solution of the problem stated
  in the beginning of section \ref{ProblemPants}. Suppose
  that the point $p(y)$ completes turns on the corresponding circle $d_r$, $d_g$
  and $d_b$ times when the argument $y$ runs around the
  'red','green' and 'blue'  boundary component of ${\cal P}(R_3)$
  respectively. We can apply the just introduced {\it construction
  1} and glue the three
  discs $\mbox{\sets U}_r$, $\mbox{\sets U}_g$, $\mbox{\sets U}_b$, to the holes of the pants.
  Essentially, we arrive at a commutative diagram:

\begin{equation}
\label{diagram}
\begin{picture}(70,50)
\put(10,40){${\cal P}(R_3)$}
\put(30,40){$\stackrel{inclusion}\longrightarrow$}
\put(50,40){$\mbox{\sets C}P^1$}

\put(15,38){\vector(1,-1){30}}
\put(18,20){$p(y)$}
\put(53,38){\vector(0,-1){30}}
\put(58,20){$\tilde{p}$ is branched}
\put(62,16){covering.}
\put(50,0){$\mbox{\sets C}P^1$}
\end{picture}
\end{equation}

\index{Riemann-Hurwitz formula}
  Applying the Riemann--Hurwitz formula for the holomorphic mapping
  $\tilde{p}$ with four ramification points (three of them are in
  the artificially glued discs and the fourth is inside the pants) we
  immediately get:

\begin{equation}
\label{RHformula}
d_r+d_g+d_b=2N,
\quad N:=\deg\tilde{p}.
\end{equation}

\subsubsection{Intersection of Circles}
\begin{lmm}
The circle $C$ does not intersect the two other circles
$\varepsilon^{\pm1}\hat{\mbox{\sets R}}$. Therefore the spectral
parameter $1<\lambda<2$ when the projective structure $p(y)$ branch
point is inside the pants.
\end{lmm}
P~r~o~o~f:~~~
  We know that the point $0$ lies in the intersection of two of
  our circles: $\varepsilon\hat{\mbox{\sets R}}$ and
  $\varepsilon^2\hat{\mbox{\sets R}}$. The total number
  $\sharp\{\tilde{p}^{-1}(0)\}$ of the pre-images of this point
  (counting the multiplicities)  is $N$  and cannot be less than
  $d_b+d_g$ -- the number of pre-images on the blue and green
  boundary components of the pants. Comparing this to
  (\ref{RHformula}) we get $d_r\ge N$ which is only possible when

\begin{equation}
\label{k1k2kRelation}
d_r=d_g+d_b=N.
\end{equation}

  Assuming that the circle $C$ intersects any of the circles
  $\varepsilon^{\pm1}\hat{\mbox{\sets R}}$ we repeat the above
  argument for the intersection point and arrive at the
  conclusion $d_b=d_r+d_g=N$ or $d_g=d_r+d_b=N$ which is incompatible with
  the already established equation (\ref{k1k2kRelation}). ~~~\bl

{\bf R~e~m~a~r~k}: In section \ref{ConeCase} we promised to show
that any  meromorphic function $p$ mapping the boundaries of the
pants to the circles (\ref{MapBoundary}) has a {\it critical
point}. Indeed, the inequalities $d_b+d_g\le N$ and $d_r\le N$
remain true whatever the branching of the structure $p$ is, while
(\ref{RHformula}) originating from the Riemann--Hurwitz formula
takes the form $d_b+d_g+d_r=2N+1$ for the {\it unbranched}
structure which leads to a contradiction. \index{projective
structure!unbranched}

\subsubsection{Image of the Pants} \label{ImPants1}

\begin{figure}[ht!]
\begin{picture}(165,55)
\put(0,0){\psfig{figure=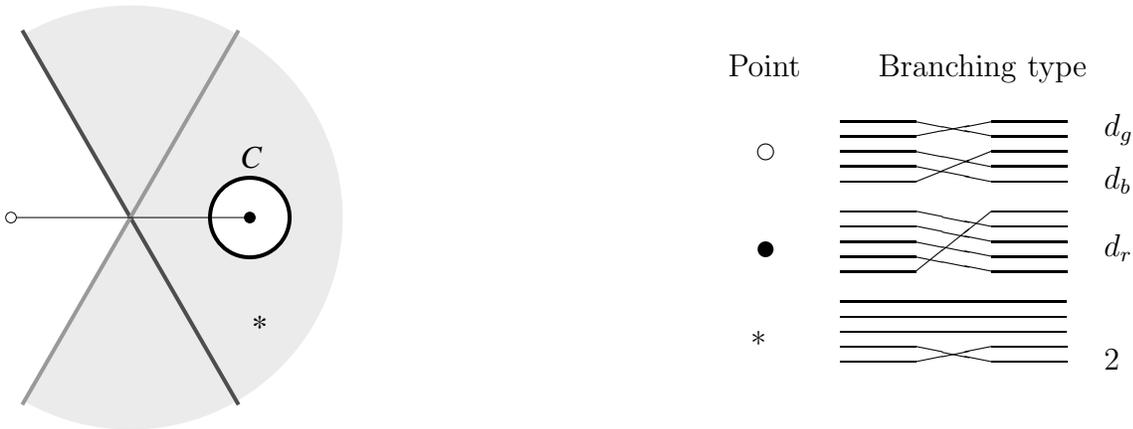}}

\put(100,0){
\begin{picture}(70,50)
\put(-5,47){Point} \put(15,47){Branching type}

\put(0,37){\circle{2}}
\multiput(10,41)(0,-2){5}{\line(1,0){10}}
\multiput(30,41)(0,-2){5}{\line(1,0){10}}
\multiput(20,37)(0,-2){2}{\line(5,-1){10}}
\put(20,41){\line(5,-1){10}}
\put(20,39){\line(5,1){10}}
\put(20,33){\line(5,2){10}}
\put(45,39){$d_g$}
\put(45,32){$d_b$}

\put(0,24){\circle*{2}}
\multiput(10,29)(0,-2){5}{\line(1,0){10}}
\multiput(30,29)(0,-2){5}{\line(1,0){10}}
\multiput(20,29)(0,-2){4}{\line(5,-1){10}}
\put(20,21){\line(5,4){10}}
\put(45,23){$d_r$}

\put(-2,10){*}
\multiput(10,17)(0,-2){3}{\line(1,0){30}}
\multiput(10,11)(0,-2){2}{\line(1,0){10}}
\multiput(30,11)(0,-2){2}{\line(1,0){10}}
\put(20,11){\line(5,-1){10}}
\put(20,9){\line(5,1){10}}
\put(45,8){$2$}

\end{picture}}

\end{picture}
\caption[] {\normalsize ({\bf a})~Shaded area is the image of
pants ~~~~~~({\bf b})~Branching type of the branch points}
\label{BelyiMap}
\end{figure}

  Let us investigate where the artificially glued discs are
  mapped to. Suppose for instance that the disc $\mbox{\sets
  U}_r$ is mapped to the exterior of the circle $C$. The point
  $0$ will be covered then at least $d_r+d_g+d_b=2N$ times which
  is impossible. The discs $\mbox{\sets U}_g$ and $\mbox{\sets
  U}_b$  are mapped to the left of the lines
  $\varepsilon{\mbox{\sets R}}$ and $\varepsilon^2{\mbox{\sets
  R}}$ respectively, otherwise points from the interior of the
  circle $C$ will be covered more that $N$ times.  The image of
  the pair of pants ${\cal P}(R_3)$ is shown on the left of the
  Fig. \ref{BelyiMap}.

  We use the ambiguity in the construction of the glueing of the disks
  to the pants and require  that the critical values of
  $\tilde{p}$ in the discs $\mbox{\sets U}_g$, $\mbox{\sets U}_b$ coincide.
  Now the branched covering $\tilde{p}$ has only three different
  branch points shown as $\bullet$, $\circ$ and $*$
  on the Fig. \ref{BelyiMap}a). The branching type at those
  three points for $d_g=2,d_b=3$, $d=N=5$ is shown on the Fig. \ref{BelyiMap}b).
  The coverings with three branch points are
  called {\it Belyi maps} and are described by certain graphs known as
  Grothendieck's {\it "Dessins d'Enfants"}. In our case the {\it
  dessin} is the lifting of the segment connecting white and
  black branch points: $\Gamma:=\tilde{p}^{-1}[\bullet,\circ]$.
\index{Grothendieck's dessins}

\subsubsection{Combinatorial Analysis of the Dessins}
\begin{figure}[ht!]
\begin{picture}(165,55)
\put(0,5){\psfig{figure=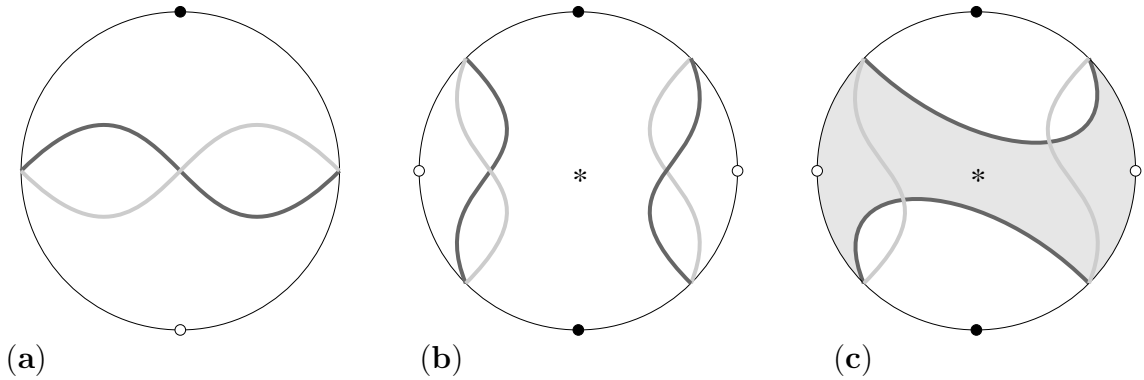}} \put(5,0){({\bf a})}
\put(60,0){({\bf b})} \put(115,0){({\bf c})}
\end{picture}
\caption[] {\normalsize ({\bf a}) Simple cell ~~ ({\bf b}) Double
cover ~~ ({\bf c}) Impossible double cover} \label{TwoCells}
\end{figure}

  There is exactly one critical point of $\tilde{p}$ over the branch point $*$.
  Hence, the complement to the graph $\Gamma$
  on the upper sphere of the diagram
  (\ref{diagram}) contains exactly one cell mapped $2-1$ to the lower sphere.
  The rest of the components of the complement are mapped $1-1$.
  Two types of cells are shown
  in figures  \ref{TwoCells} a) and  b), the
  lifting of the red circle is not shown to simplify the
  pictures. The branch point $*$ should lie in the intersection of
  the two annuli $\alpha$ and $\overline{\alpha}$, otherwise
  the discs $\mbox{\sets U}_g$, $\mbox{\sets U}_b$ glued to different boundary components of
  our pants will intersect: the hypothetical case when the
  branch point of $p(y)$ belongs to one annulus but does not
  belong to the other is shown in Fig. \ref{TwoCells} c).

\begin{figure}
\psfig{figure=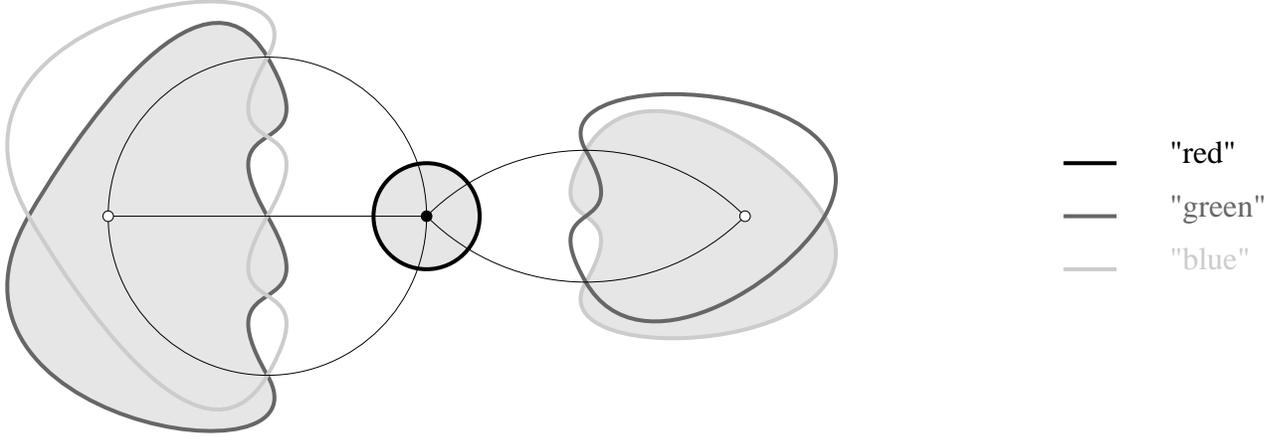} \caption[] {\normalsize Dessin for
$d_g=3,d_b=2$; the pre-image of the branch point $*$ is at the
infinity} \label{Assembly1}
\end{figure}

  The cells from Fig. \ref{TwoCells} a), b) may be assembled in
  a unique way shown in Fig.~\ref{Assembly1}. The pants are
  colored in  white, three artificially sewed discs are
  shaded. Essentially this picture shows us how to sew together
  the patches bounded by our three circles $C$,
  $\varepsilon^{\pm1}\hat{\mbox{\sets R}}$ to get the pants
  conformally equivalent to ${\cal P}(R_3)$. As a result of the
  surgery procedure we obtain
  the pants ${\cal P}_1(\lambda, h_1, h_2|d_g-1,d_b-1)$.

\subsection{Simple branch points on the boundary of the pants}\label{OnBorder}
Our strategy remains the same: to fill in the holes in the pants and to convert $p(y)$
into a branched covering with a simple type of branching.

\subsubsection{Construction 2}\label{Construct2}
  Let again $p(y)$ be a holomorphic mapping of a bounded Riemann
  surface $\cal M$ to the sphere with the selected boundary
  component $(\partial {\cal M})_*$ being mapped to the boundary
  of the unit disc {\sets U}. Now the mapping $p(y)$  has two
  simple critical points on the selected boundary component (the
  case of coinciding critical values is not excluded). Those two
  points divide the oval $(\partial {\cal M})_*$ into two
  segments: $(\partial {\cal M})_*^+$ and $(\partial {\cal M})_*^-
  $. Let the increment of $\arg ~p(y)$ on the segment $(\partial
  {\cal M})_*^+$ be $2\pi d^+-\phi$, $0<\phi\le2\pi$, and the
  decrement on the segment $(\partial {\cal M})_*^-$ be $2\pi d^--
  \phi$,  the point $y$ moves around the selected
  oval in the positive direction and $d^\pm$ are positive
  integers. We are going to fill in the hole in the Riemann
  surface $\cal M$ with two copies  of the unitary disc (\ref{Udisc}): $\mbox{\sets U}^+$ and
  $\mbox{\sets U}^-$.\index{winding number}

We define the mapping from the disjoint union ${\cal M}\cup$  $\mbox{\sets U}^+\cup$  $\mbox{\sets U}^-$
to the sphere:
\begin{equation}
\tilde{p}(y):=\left\{
\begin{array}{ll}
p(y),& y\in {\cal M},\\
L^-(y^{d^-}),& y\in \mbox{\sets U}^-,\\
L^+(y^{-d^+}),& y\in \mbox{\sets U}^+,\\
\end{array}
\right.
\end{equation}
where $L^\pm(\cdot)$ are the (at the moment arbitrary) linear
fractional mappings
keeping the unitary disc (\ref{Udisc}) invariant. The choice of $L^\pm(\cdot)$
will be specified later to simplify the combinatorial analysis.

\begin{figure}
\begin{picture}(165,60)
\thicklines
\put(30,25){\oval(40,40)}
\thinlines
\put(30,5){\line(0,1){40}}
\put(29,4){$*$}
\put(29,44){$*$}

\put(0,25){({\bf a})} \put(20,25){$\mbox{\sets U}^+$}
\put(35,25){$\mbox{\sets U}^-$} \put(55,25){$\cal M$}

\put(10,48){$(\partial{\cal M})_*$}

\put(100,25){({\bf b})}
\put(120,0){\psfig{figure=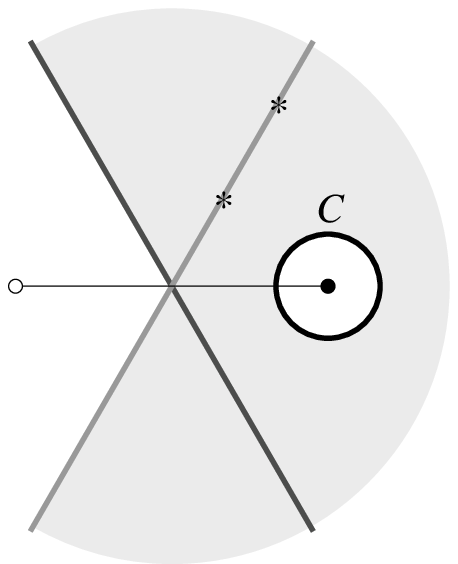}}

\end{picture}

\caption[] {\normalsize ({\bf a}) Filling in the hole bounded by
$(\partial{\cal M})_*$ \hfill ({\bf b}) The shaded area is the image of
${\cal P}(R_3)$.} \label{GlueTwoDiscs}
\end{figure}

  Identifying the points $y$ with the same value of
  $\tilde{p}(y)$ we glue the segments $(\partial {\cal M}*)^\pm$
  of the selected boundary oval of $\cal M$ to the portions of
  the boundaries of the discs  $\mbox{\sets U}^\pm$ respectively.
  The remaining parts of the boundaries of $\mbox{\sets U}^\pm$
  are glued to each other as shown in Fig.
  \ref{GlueTwoDiscs}a).

\begin{figure}
\begin{picture}(165,50)
\thicklines
\put(30,40){\oval(40,10)}
\put(5,40){\oval(10,60)[br]}
\put(55,40){\oval(10,60)[bl]}
\thinlines
\put(30,10){\oval(50,13)[b]}
\put(20,35){\vector(-1,0){3}}
\put(40,45){\vector(1,0){3}}
\put(30,34){$*$}
\put(30,44){$*$}
\put(30,10){$\cal M$}
\put(10,48){$(\partial{\cal M})_*$}
\put(75,23){\vector(1,0){20}}
\put(80,25){$p(y)$}

\put(120,5){\psfig{figure=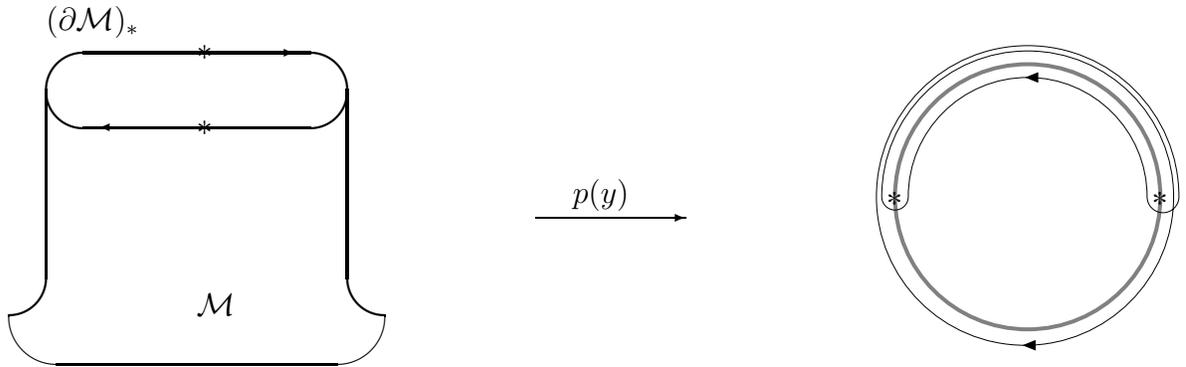}}

\end{picture}

\caption[] {\normalsize Mapping of the boundary component
$(\partial{\cal M})_*$ with two simple branch points $*$ on it and
winding indices $d^+=1, d^-=2$.} \label{BdryCritPts}
\end{figure}

\subsubsection{Branched Covering of a Sphere}

  At the moment we do not know which of the three boundary ovals
  of the pants ${\cal P}(R_3)$ contains the critical points of
  $p(y)$. Therefore we introduce the 'nicknames' $\{1,2,3\}$
  for the set of colors $\{r,g,b\}$ so that the critical points
  will be on the oval number 3. The usage of {\it construction 2} from
  section  \ref{Construct2} allows us to glue two discs
  $\mbox{\sets U}^\pm_3$ to the latter boundary. The usage of
  {\it construction 1} from section  \ref{Construct1} fills in the
  remaining two holes with two discs $\mbox{\sets U}_1$ and
  $\mbox{\sets U}_2$. Positive integers arising in those
  constructions are denoted by $d_3^\pm$, $d_1$, $d_2$ respectively.

Again, we split the mapping $p(y)$ from the pants to the sphere as
in the diagram (\ref{diagram}): $p=\tilde{p}\circ{~inclusion}$ with
the branched covering $\tilde{p}$. The latter mapping has six
critical points: two
  simple ones inherited from the pants and four at the centers of the
  artificially glued discs and multiplicities $d_3^\pm-1$, $d_1-
  1$, $d_2-1$ respectively. The Riemann--Hurwitz formula for this covering gives
 \begin{equation}
 \label{RH2}
  d_1+d_2+d_3^++d_3^-=2N,
  \qquad N:=\deg\tilde{p}.
 \end{equation}

\begin{lmm}
The images of the ovals with numbers 1 and 2 do not intersect.
\end{lmm}

P~r~o~o~f. Suppose the opposite is true and a point $Pt$ lies in the
intersection of the images of the first two ovals. Then
$N\ge\sharp\tilde{p}^{-1}(Pt)\ge d_1+d_2$. On the image of the third
oval there is a point (e.g. in the right side of  Fig.
\ref{BdryCritPts} this is a point $i$) with $d_3^++d_3^-\le N$
pre-images. Comparing the last two inequalities to (\ref{RH2}) we
get the equalities
$$
d_1+d_2=d_3^++d_3^-=N
$$
and $Pt$ is covered at least $d_1+d_2+\min(d_3^+,d_3^-)>N$ times. ~~\bl

{\bf Corollary.}
Two circles $\varepsilon^{\pm1}\hat{\mbox{\sets R}}$ intersect, therefore the
critical points of $p(y)$
lie either on the blue or on the green boundary of pants.
Moreover, the circle $C$ --  the image of the red boundary oval --
does not intersect the two mentioned circles which may only happen when
$\mu\in(\frac12,1)$, or equivalently $\lambda\in(1,2)$.

{\bf Convention:} We assume that both critical points of $p$ lie on the
blue oval. The remaining case when they belong to the green oval is
absolutely analogous to the case we consider. Now the notations
$\mbox{\sets U}_b^\pm$, $\mbox{\sets U}_r$, $\mbox{\sets U}_g$,
$d_b^\pm$, $d_r$, $d_g$ have the obvious meaning.

\subsubsection{The Image of the Pants} \label{ImPants2}
Let us show that the the image of the pants remains the same as in
section \ref{ImPants1}.

\begin{lmm}
The image $p({\cal P}(R_3))$ of the pants
lies in the intersection of annuli $\alpha$ and $\bar{\alpha}$
  -- see Fig. \ref{GlueTwoDiscs}b)
\end{lmm}

P~r~o~o~f. We refer to the four sectors:  {\sets C}$\setminus
\varepsilon^{\pm1}${\sets R} as to 'top', 'down', 'left' and
'right'. It is a matter of notation to say that the disc
$\mbox{\sets U}_b^+$ is mapped to the 'top' and 'left' sectors while
the disc $\mbox{\sets U}_b^-$ is mapped to the 'down' and 'right'
sectors.

  The disc $\mbox{\sets U}_g$ covers either the 'top' or the 'left'
  sector and both are covered by the disc $\mbox{\sets U}_b^+$.
  Therefore, $d_g+d_b^+\le N$. In a similar way we get
  $d_r+d_b^-\le N$. The obtained inequalities and the Riemann--Hurwitz formula
  (\ref{RH2}) -- which in our notations becomes
  $d_r+d_g+d_b^++d_b^-=2N$ -- give us
$$
d_r+d_b^-=d_g+d_b^+=N.
$$

If the disc $\mbox{\sets U}_r$ is mapped to the exterior of the circle $C$, then
either 'left' or 'top' sector is covered $d_r+d_g+d_b^+>N$ times.
If the disc $\mbox{\sets U}_g$ is mapped to the right of the line $\varepsilon${\sets R}, then
 the interior of the circle $C$ is  covered $d_r+d_g+d_b^->N$ times.

We see that the 'left' sector and the interior of the circle $C$ are covered by the
artificially inserted discs only. ~~~\bl

{\bf Corollaries} ~~1) Both critical values of $p(y)$ lie on the ray
$-\varepsilon^2(0,\infty)$. ~~2) The integer $d_b^-$ is equal to $1$, since the point $0$ is covered at least
$d_g+d_b^++d_b^--1\le N$ times.\\[2mm]

Let us recall that the constructions of attaching discs to the pants
allow us to move branch point (= the critical value of $\tilde{p}(y)$ in
the inserted disc) within the appropriate circle. In particular, the
critical values of $\tilde{p}(y)$ in the discs $\mbox{\sets U}_g$,
$\mbox{\sets U}_b^+$ may be placed to the same point in the 'left'
sector, say to $p=-1$ (point $\circ$ in Fig.
\ref{GlueTwoDiscs}b) while the critical values in the discs
$\mbox{\sets U}_r$, $\mbox{\sets U}_b^-$ may be placed to the same
point inside $C$, say to $p=1$ ( point $\bullet$ in Fig.
\ref{GlueTwoDiscs}b).  Now we lift the segment $[\circ, \bullet]$
connecting the branch points to the upper sphere of the diagram
(\ref{diagram}) and analyze the arising graph
$\Gamma:=\tilde{p}^{-1}([\circ, \bullet])$.

\subsubsection{Combinatorial Analysis of the Graph}
The restriction of $\tilde{p}$ to every component $F$ of the
compliment $\hat{\mbox{\sets C}}\setminus\Gamma$ to the graph is
naturally continued to the branched coverings over the
disc\footnote{{\it Closure} here has the same meaning as in the
formula (\ref{R3Pants}) } $Closure(\hat{\mbox{\sets
C}}\setminus[\circ, \bullet])$. We can list all flat surfaces $F$
covering a disc with the branching number $B\le2$. To this end we
use the Riemann--Hurwitz formula for the branched coverings of the
bordered surfaces: \index{Riemann-Hurwitz formula}
$$
2+B=\sharp\{\partial F\}+\deg\tilde{p}|F
$$
which relates $B$ -- the total branching number of $\tilde{p}$ in
the selected flat surface $F$ covering a disc; $\sharp\{\partial
F\}$ -- the number of its boundary components and $\deg\tilde{p}|F$
-- the degree of the restriction of the covering $\tilde{p}$ to the
component $F$. Taking into account that $\sharp\{\partial
F\}\le\deg\tilde{p}|F$ we obtain the list shown in
Tab.~\ref{DiscCover}.

\begin{table}
\caption{Flat surfaces $F$ covering the disc with the branching
number $B\le2$} \centering \label{DiscCover}
\begin{tabular}{c|c|c|c}
number of sheets & $B$ &
surface $F$ & picture \\
\hline
1&0&disc&Fig. \ref{TwoCells}({\bf a})\\
2&1&disc& Fig. \ref{ComplimentSurface}({\bf a})\\
3&2&disc& Fig. \ref{ComplimentSurface}({\bf b})\\
2&2&annulus& Fig. \ref{ComplimentSurface}({\bf c})\\
\end{tabular}
\end{table}

\begin{figure}
\begin{picture}(165,55)
\put(0,5){\psfig{figure=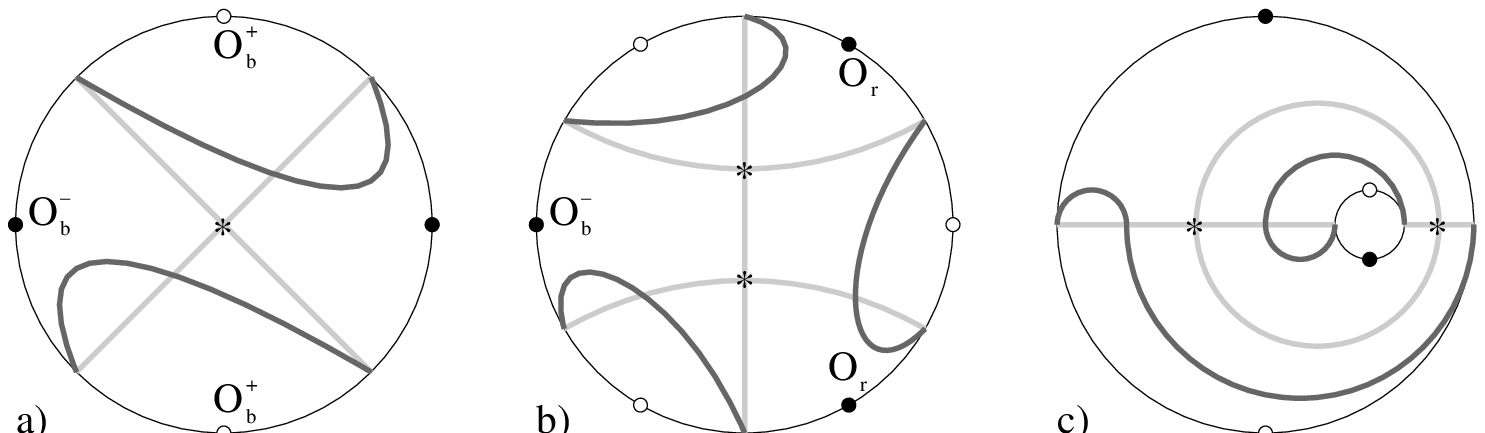}} \put(5,0){({\bf a})}
\put(60,0){({\bf b})} \put(115,0){({\bf c})}
\end{picture}
\caption[] {\normalsize Flat surfaces $F$ covering the disc with the
branching numbers $B=1,2$.} \label{ComplimentSurface}
\end{figure}

\begin{figure}
\centerline{\psfig{figure=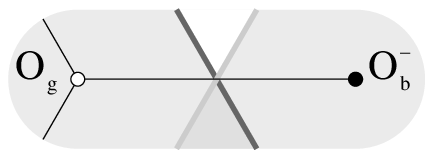}} \caption[]
{\normalsize If $O_g$ and $O_b^-$ were neighbors, the discs
$\mbox{\sets U}_b^-$ and $\mbox{\sets U}_g$ would intersect.}
\label{ImposibleNeighbour}
\end{figure}

  The combinatorics of the green and blue circles lifted to the listed covering surfaces $F$
  is shown in the Fig. \ref{TwoCells}a) and
  Fig.\ref{ComplimentSurface}a-c). Let us denote the centers of the four
  artificially glued discs $\mbox{\sets U}_r$, $\mbox{\sets
  U}_g$, $\mbox{\sets U}_b^+$ $\mbox{\sets U}_b^-$ as
  respectively $O_r$ (black vertex of graph $\Gamma$ with valency
  $d_r$), $O_g$, $O_b^+$ (white vertexes with valencies $d_g$,
  $d_b^+$) and $O_b^-$ (dangling black vertex). Their
mutual positions in the graph $\Gamma$ are subject to the following restriction:

\begin{lmm}
The vertices $O_g$ and $O_b^-$ are not neighbors in $\Gamma$.
\end{lmm}
P~r~o~o~f: The disjoint discs $\mbox{\sets U}_b^-$ and $\mbox{\sets U}_g$ of the upper
sphere in the diagram (\ref{diagram}) would intersect otherwise -- see Fig. \ref{ImposibleNeighbour}. ~~\bl

{\bf Corollaries} 1)
The vertices on the border of the triply covering disc $F$
-- see Fig. \ref{ComplimentSurface}b) --
appear in the following order: $O_g$, $O_r$, $O_b^+$, $O_b^-$, $O_b^+$, $O_r$.
They may be uniquely ascribed to the vertices in the picture
after the  observation:
{\it the blue line divides the vicinity of any critical point $*$ into
four quadrants, two of which belong to the pair of pants, one belongs to the disc
$\mbox{\sets U}_b^-$, and the rest is contained in the disc $\mbox{\sets U}_b^+$}.

2) The complement to the graph $\Gamma$ cannot contain two doubly
covering discs $F$. Indeed, the point $O_b^-$ lies on the boundary
of one of those discs. Both neighboring vertices on the boundary
of the disc $F$ should be $O_b^+$ according to the lemma. But this
contradicts the above {\it observation}: two quadrants of this
covering disc belong to $\mbox{\sets U}_b^+$ -- see Fig.
\ref{ComplimentSurface}a).

\subsubsection{Assembly Scheme}
\begin{figure}
\begin{picture}(165,55)
\put(-2,5){\psfig{figure=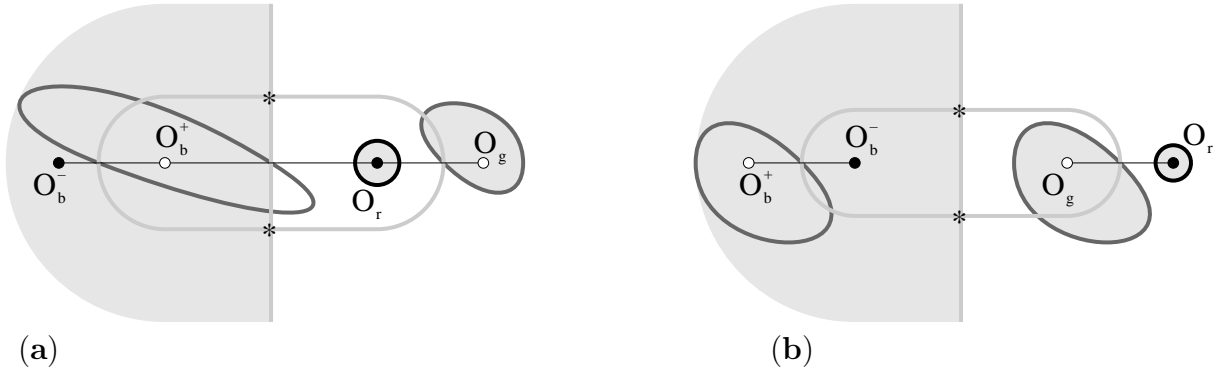}} \put(0,0){({\bf a})}
\put(100,0){({\bf b})}
\end{picture}
\caption[] {\normalsize Graph $\Gamma$ for the basic mappings with
$d_g=1$, $d_b^+=d_r=2$, $N=3$ ({\bf a}) and $d_g=d_b^+=d_r=1$, $N=2$
({\bf b}). Artificially inserted discs are shaded.}
\label{Assembly2}
\end{figure}

We see that there remain only two possibilities for the complement
to the graph $\Gamma$. It consists either of ({\bf a}) one disc
mapped 3-1 and $N-3$ {\it simple} cells mapped 1-1 or ({\bf b}) an
annulus mapped 2-1 and $N-2$ simple discs mapped 1-1. The graphs
$\Gamma$ with compliment containing no simple cells are shown in
Fig. \ref{Assembly2}. They correspond to the pants ${\cal
P}_2(\dots|0,1)$ ({\bf a}) and ${\cal P}_2(\dots|0,0)$ ({\bf b}).
The graphs with simple cells in the complement are obtained from
those two basic pictures as a result of the surgery. We cut the
graph along the the edge $O_rO_g$ and insert $d_g-1$ simple discs
in the slot as in Fig. \ref{Assembly1}. The graph on the left side
of the Fig. \ref{Assembly2} admits another surgery: we cut the
graph along the edge $O_rO_b^+$ and sew in $d_b^+-2$ patches
shown in Fig. \ref{TwoCells}a) in the slot. The arising graph
corresponds to the pair of pants ${\cal P}_2(\dots,
d_g-1,d_b^+-1)$.

\subsection{Remaining cases}
  If the branch points of the projective structure $p:=p^+$ belong
  to the green oval of the pants we arrive at the pair of pants ${\cal P}_s$
  of fashion $s=3$. Finally, when the branch points merge the limit
  variant of construction 2 may be applied for the analysis and we
  arrive at the pants of intermediate types $s=12,~13$.

\section{Conclusion}
A similar analysis based on the geometry and combinatorics may be
applied to obtain the representations of the solutions of the PS-3
integral equation in all the omitted cases. Much of the techniques
used may be helpful for the study of other integral equations with
low degree rational kernels.


\vspace{5mm}
\parbox{8cm}
{\small\tt
119991 Russia, Moscow GSP-1,\\
ul. Gubkina 8,\\
Institute for Numerical Mathematics,\\
Russian Academy of Sciences}

\end{document}